\documentclass[12pt]{article}

\usepackage{amsmath}
\usepackage{amsfonts}
\usepackage{amsthm}
\usepackage{epsfig}
\usepackage{cite}
\usepackage{color} 

\newcommand \comment[1]{}                       
\renewcommand \comment[1]{\emph{#1}}            

\newtheorem{theorem}{Theorem}[section]
\newtheorem{proposition}[theorem]{Proposition}
\newtheorem{corollary}[theorem]{Corollary}
\newtheorem{lemma}[theorem]{Lemma}
\theoremstyle{definition}
\newtheorem{definition}[theorem]{Definition}

\newcommand{\Remark}{\textbf{Remark}}
\newcommand{\Remarks}{\textbf{Remarks}}

\newcommand{\MVAR}[1]{{[#1]\;}}

\newcommand{\dunion}
{\mbox{\hbox{\hskip6pt$\cdot$\hskip-5.50pt$\cup$\hskip2pt}}}
\newcommand{\scomp}[1]{\ensuremath{\overline{#1}}}

\newcommand{\scma}{\ensuremath{\ }}

\newcommand{\epsilonfun}{\ensuremath{\epsilon}}

\newcommand{\Is}{\ensuremath{\iota}}
\newcommand{\Vs}{\ensuremath{\upsilon}}

\newcommand{\Reals}{\ensuremath{\mathbb{R}}}
\newcommand{\FieldK}{\ensuremath{K}}
\newcommand{\Perms}{\ensuremath{\mathfrak{S}}}
\newcommand{\rank}{{\rho}}
\newcommand{\Rank}{{\rho}}
\newcommand{\Card}[1]{\ensuremath{{\left|#1\right|}}}
\newcommand{\ext}[1]{\ensuremath{\mathbf{#1}}}
\newcommand{\extvee}{\;\;}
\newcommand{\Plucker}{Pl\"{u}cker\ }

\makeatletter
\newfont{\footsc}{cmcsc10 at 8truept}
\newfont{\footbf}{cmbx10 at 8truept}
\newfont{\footrm}{cmr10 at 10truept}
\makeatother
\title{Ported Tutte Functions of Extensors and Oriented Matroids}

\author{Seth Chaiken
\thanks{Some of this research was 
performed during my 2001 academic year 
Sabbatical at the University at Albany.}\\
\small Computer Science Department\\[-0.8ex]
\small University at Albany, State Univ. of New York, USA\\[-0.8ex]
\small \texttt{sdc@cs.albany.edu}
}

\date{August, 2 2006}

\begin{document}

\maketitle

\newpage
\begin{abstract}
The Tutte equations are ported (or set-pointed) when the equations
$F(\mathcal{N})$ $=$ $g_eF(\mathcal{N}/e)$ $+$ 
$r_eF(\mathcal{N}\setminus e)$ are 
omitted for elements $e$ in a distinguished set
called ports.
The solutions $F$, called ported Tutte
functions, can distinguish different orientations of the same matroid. 
A ported extensor with ground set is a (fully) decomposable element
in the exterior algebra (of antisymmetric tensors) over a vector space with
a given basis, called the ground set, containing a distinguished subset 
called ports.  A ported extensor is one way to present a linearly 
representable ported matroid or oriented matroid.  There are extensor operations
corresponding to oriented matroid dualization, and to deletions and
contractions.

We define a ported extensor function by means of dualization, 
port element renaming, exterior multiplication, 
and then contraction of all non-port elements.  The main result is
that this function  satisfies 
a sign-corrected variant of the Tutte equations in which deletion
and contraction are extensor operations, and addition and 
the anticommutative multiplication belong to an exterior algebra
rather than to a commutative ring.

For extensors representing unimodular, i.e., regular
matroids with an empty port set, our function reduces to the basis
generating function; and then, for graphs, to the Laplacian (or
Kirchhoff) determinant.  On graphs with port edges, the function
value, as an extensor, signifies the space of solutions to Kirchhoff's
and Ohm's electricity equations after projection to the voltage and
current variables associated to the ports.  
In
particular, the Laplacian matrix with the identity matrix appended 
presents one example of our extensor function's value.  
Combinatorial interpretation
of various determinants (the \Plucker coordinates) generalize 
the matrix tree theorem and forest enumeration expressions for 
electrical resistance.  

We also demonstrate how the corank-nullity polynomial, basis expansions
with activities, and a geometric lattice expansion generalize to 
ported Tutte functions of oriented matroids.  The ported 
Tutte functions are parametrized, which raises 
the problem of how to generalize known characterizations
of parametrized non-ported Tutte functions.
\end{abstract}

\newpage

\tableofcontents

\section{Introduction}

A ported unoriented or oriented matroid 
$\mathcal{N}=\mathcal{N}(P,E)$ 
has its ground set 
$S(\mathcal{N})=P\dunion E$ given with a distinguished
subset $P$ of elements which we call \textbf{ports}; $P\cap E=\emptyset$.
The following definition 
combines the idea of Tutte invariants of set-pointed matroids
studied by Las Vergnas
\cite{MR0419272,TPMorphMat80,SetPointedLV}
with the idea of the parametrized Tutte equations and functions 
studied by 
Zaslavsky\cite{MR93a:05047} (who calls our Tutte functions ``strong''), 
Bollobas and 
Riordan\cite{BollobasRiordanTuttePolyColored} and
Ellis-Monaghan and Traldi \cite{Ellis-Monaghan-Traldi}.
Let two parameters $g_e$ and $r_e$ be given for each $e\in E$.

\begin{definition}
\label{firstTutteFunction}
A function $F$ is a 
\textbf{(ported and parametrized) Tutte function}
if the domain of  $F$ is 
a minor closed class
of ported unoriented or oriented matroids and
$F$ satisfies the following
\textbf{(ported and parametrized) Tutte Equations}
for each $\mathcal{N}$ in the class:

When $e\in E$ is a non-separating element, i.e.,
$e$ is neither a port nor loop nor a coloop (i.e., isthmus):
\begin{equation}\label{TGplus}
	F(\mathcal{N}) = g_eF(\mathcal{N}/ e) + r_eF(\mathcal{N}\setminus e).
\end{equation}
When $\mathcal{N}=\mathcal{N}_1\oplus\mathcal{N}_2$
with ground sets $S(\mathcal{N}_1)\cap S(\mathcal{N}_2)=\emptyset$:
\begin{equation}\label{TGtimes}
	F(\mathcal{N}_1\oplus\mathcal{N}_2) = F(\mathcal{N}_1)F(\mathcal{N}_2).
\end{equation}

\end{definition}

In the following, a Tutte function shall be assumed to be 
both ported and parametrized as above unless otherwise indicated.  The letters
$g_e$ and $r_e$ which we use for parameters, and the terminology
``set of ports'' for the distinguished ground set elements,
are suggested by the electrical network theory 
application which motivates our research.  Briefly, $g_e$ (for conductance
or admittance)
and $r_e$ 
(for resistance or impedance\footnote{Electrical engineers
customarily use the words admittance and impedance for complex values.}) 
appear in the homogeneous 
expression $r_e:g_e$ of a resistance value (i.e., the resistance value
is either $r_e/g_e$ ohms if $g_e\neq 0$, or 
infinity if $g_e=0$ and $r_e\neq 0$.)
A port
designates one interface 
(a pair of terminal vertices and two variables, 
one for current and one for voltage) between the physical electrical network
and an external environment.  (Think of an ordinary 2-pin electrical
plug or receptacle.)  
Other physical applications of parametrized Tutte functions of graphs are 
surveyed in \cite{Sokal}.  
Details and references for the electrical application are presented
but the logical background, results and proofs 
(\textsection \ref{Preliminaries}-\ref{Activity})
are independent of these details.

Our main topic is a generalization of one construction of a 
linear subspace from the graph (with port edges and parameters)
that occurs in the application.  The constructed
subspace consists of the space of network solutions projected to
the port variables.  We generalize by replacing the graph's coboundary
(cocycle) space by a finite dimensional linear subspace.  The underlying
space has a particular ground set basis, whose elements
generalize graph edges, and this ground set has the distinguished 
subset of port elements.

Our formulation with ports subsumes 
(see \textsection\ref{Laplacian}) the classical 
equation $L(g,1)\phi = J$ on the node voltages $\{\phi_i\}$ and 
external currents flowing into nodes $\{J_i\}$.  Here, $L(g,1)$ denotes the
edge-weighted Laplacian matrix; the edge weights $g_e$ are the conductances
(or admittances) of edges and we take each $r_e=1$.  Assuming the graph is
connected, each principal cofactor of $L(g,1)$ equals (according to the
famous Matrix Tree Theorem) the weighted spanning
tree polynomial
\[
T(g)=\sum_{
\substack{T\subseteq E\\
T\text{ a spanning tree}
         }
          }
\prod_{e\in T}g_e
=
\sum_{T}g_T.
\]

Choe, Oxley, Sokal and Wagner\cite{ChoeOSW} prove that 
this and similar basis generating polynomials for some but not
all other matroids have the 
\textbf{half-plane property}:
For complex
values $g_e$, if $\Re g_e > 0$ for all $e\in E$ then $T(g)\neq 0$.
The half-plane property 
characterizes those electrical networks that 
cannot generate energy.
These authors show that the monomials of every non-zero homogeneous
multiaffine polynomial with the half-plane property enumerate the bases of
a matroid.  After noting that the half-plane property is a strengthening
of the property that all the coefficients have the same phase, they raise
the open question of whether, for every matroid, coefficient values
exist so the basis enumerating polynomial with these coefficients
has the half-plane property.  It is noted that such coefficients can
always be assumed to be positive.

We treat polynomials 
(including $T(g)$ when $r_e=1$) whose ratios,
for graphs, signify externally observable numerical characteristics of 
an electrical network that pertain to several port current or 
voltage values together.
These quantities include the 
``transpedances'' used by
Brooks, Smith, Stone and Tutte\cite{BSST} in their work on 
finding dissections of squares into squares by means of 
solving electrical network equations.  Indeed, these authors'
presentation of the Matrix Tree Theorem is in the context
of combinatorial proofs of solution properties.  Our treatment
describes the coefficient signs in terms of oriented matroids.
Some of the polynomials have the half-plane property and
some clearly do not because they have terms with different signs.
We therefore hope that our work may contribute insight into
the questions posed in \cite{ChoeOSW}, as well as more recent questions
raised by Wagner\cite{NegCorr}
pertaining to inequalities among differences between 
polynomials that enumerate certain trees or forests.  This work
combines Rayleigh's inequality (motivated by the physics of 
electrical networks with positive edge conductances, see also
\cite{ChoeRayl}) and 
analysis of correlations where ratios of edge parameters
signify probabilities.

It may help some readers to
know that the generalization of the coboundary space
and the parameters $r_e$ and $g_e$, which
together determine our constructed solution subspace, are specified
separately---This separation corresponds to the distinction between two 
kinds of constraints (exact or geometric versus approximate)
that is explained and modeled with matroid
theory by Murota\cite{Murota}.   

Our work
distinguishes the polynomial determinants where all terms 
(in the parameters) have the same sign from those
where differing signs occur.  So, a non-zero determinant with
terms of differing signs might vanish for critical combinations
of parameter values; whereas a determinant with terms all the same
sign will never vanish.  In Murota's model, this distinction 
would depend on the exact or geometric constraints; we express the 
distinction in terms of its oriented matroid properties.  

Some of our previous work \cite{sdcOMP} 
applied oriented matroids to distinguish the case of
a vanishing determinant for critical parameter value 
combinations.  The approach did not apply Tutte function
theory.  Instead, we investigated the property of
a pair of oriented matroids with a common ground set
have a common non-zero covector.  The electrical network
applications treated were more general than in the 
present paper.  We began with a dual oriented 
matroid pair (graphic and cographic), 
deleted or contracted certain elements in each, and then
evaluated the above common covector property in the resulting
oriented matroid pair with a common ground set.
See \cite{WelshKayibiLinking}
for Tutte theory developed for paired matroids.  We are 
currently investigating how to generalize paired matroids
by adding port elements and relate the theories 
of the present paper with \cite{sdcOMP} and
\cite{WelshKayibiLinking}.

\subsection{Exterior Algebra}

Exterior algebra is used to represent and operate on linear subspaces.
A (fully or completely) decomposable element in an exterior algebra is
either a field element or the exterior product of vectors.  For the sake of
brevity, such decomposable will be called \textbf{extensors}.  An extensor
is \textbf{ported} when an underlying space is given with a ground set basis
and a distinguished port subset.
We remind the reader that each non-zero extensor
corresponds to the unique subspace whose bases are the sets of vectors whose
exterior product equals the given extensor, up to a non-zero field 
element multiple (see Theorem \ref{ExtTheory}.)

For the electrical network application, a graph's coboundary space
is represented by extensor $\ext{N}$.  In our theory,
$\ext{N}$ would represent an arbitrary linear subspace.
We will define some functions of $\ext{N}$ or of certain
matroids with values within an exterior algebra.  
Each function value represents
the solution subspace for a generalization of the electrical network 
problem.
The domain of the first such function 
(see Definition \ref{MENDef})
consists of ported extensors.  
The second function applies to the class of ported unimodular
(i.e., regular) matroids.  The second function is defined by
specializing the first function to the extensors that represent
unimodular matroids.  See Definition \ref{DefMUOM}.
When the latter function is specialized to graphic oriented matroids with no
port elements, it reduces to the Laplacian determinant which, according to
the famous Matrix Tree Theorem, equals the spanning tree polynomial
discussed above.

The additive Tutte equation pertains to a graph or matroid and its
two minors obtained by deletion and contraction.  We will define algebraic 
operations on $\ext{N}$.
Their values are
denoted by $\ext{N}/e$ and $\ext{N}\setminus e$ for each $e\in E$.
We also adapt the Hodge star operation to define an operation corresponding
to oriented matroid dualization.
Our main results (Theorems \ref{TheoremMain} and \ref{EquationsExtFunMatroids})
are that these two functions obey equations similar to the 
Tutte equations above, taken in an exterior algebra rather than in
a commutative ring.  The dualization operation is used in defining
these functions and in proving the main results.
We hope the reader will bear with us in
using the terms Tutte equations and functions in this context 
before the precise equations can be presented.
If not, one can skip to the definitions and theorems in 
\textsection \ref{Preliminaries} 
and \textsection \ref{Extensors} without loss of logical continuity.
Sign factors are required in our Tutte equation variants
to accommodate the anticommutative multiplication.

It is well-known and easy to verify
that the spanning tree count is a Tutte invariant of graphs.  
Our results 
further elucidate the relationships between Tutte functions,
the Matrix Tree Theorem theorem and 
enumeration methods
for resistive electrical network solutions
pioneered by Kirchhoff\cite{Kirchhoff} and 
Maxwell\cite{MaxwellsFramesPaper}.  
These methods were introduced into combinatorial theory
by Brooks, Smith, Stone and Tutte\cite{BSST,TutteBook} who attributed
them to Kirchhoff, see \textsection \ref{Laplacian}.  
They continue to be applied
within some electrical engineering computer aided
design tools\cite{someEEBooks}.
Our two points of departure from \cite{BSST} are to 
replace analysis in terms of graph vertices and incidences
by analysis
of functions on the graph edges, and then to
express the relevant equations 
(\textsection \ref{ElecNetEquatSec}) in exterior algebra.

The extensor functions that we study are generalizations of the determinant
of the (reduced) edge-weighted Laplacian matrix of a graph.  
The homogeneous form of this determinant, according to the Matrix Tree Theorem,
equals
\begin{equation}
\label{BasisEnumPoly}
r_S\det(NDN^t) = r_S\sum_{X\subseteq S}\frac{g_X}{r_X}N[X]^2
= r_S\sum_{B}\frac{g_B}{r_B}
= \sum_{B} g_B r_{\scomp{B}};\;\;D=\mbox{diag}(\frac{g_e}{r_e}, e\in E)
\end{equation}
where $N$ is the reduced signed incidence matrix of the graph and the sum
is over all spanning spanning trees $B$.
We remind 
the reader that every extensor can be represented by a sequence of determinants
which is called its \textbf{\Plucker coordinates}.
Of course, the maximal forests
are the bases of the graphic matroid.

In summary, our theorems about extensor functions satisfying
sign-corrected Tutte equations generalize the graphic matroid 
case of the easy-to-prove fact that the basis generating function
is a (non-ported) Tutte function defined for all matroids.  

\subsection{New Tutte-Like Invariants}

However, somewhat deeper theory involving
the Laplacian, matroid invariants and exterior algebra is 
involved.  First, consider 
Tutte invariants of matroids, such as the number of bases. 
The universal Tutte invariant, the well-known Tutte polynomial,
is defined for all matroids; hence no Tutte invariant depends
on the orientation of an oriented matroid.  Let us extend the
definition as follows:
A \textbf{ported Tutte invariant} of ported oriented matroids
is a Tutte function $f$ where all the parameters 
$r_e=g_e=1$ 
that is invariant
under oriented matroid isomorphisms $f$ that 
preserve each port element.
Specifically, $f(p)=f^{-1}(p)=p$ for every port $p$ in the
domain or range of $f$.  
Las Vergnas 
\cite{MR0419272,TPMorphMat80,SetPointedLV}
developed the analogous ideas for ported matroids (without orientation),
and called them ``set-pointed.''  This work on invariants of 
matroid morphisms, i.e., strong maps, applies the
universal ported invariant called the \textbf{big Tutte polynomial}.  
We had applied this idea to study the behavior
of such polynomials under the restricted matroid union 
operation\cite{sdcPorted} where the composed matroids can only have port elements in common, 
which is one generalization of matroid series connection.  

It is
easy to develop the corresponding universal ported \emph{oriented} matroid invariant.
Since the latter has one variable for each connected oriented 
matroid whose ground set
contains port elements exclusively, we conclude that some ported Tutte
invariants of oriented matroids \emph{do} distinguish different 
orientations of the same matroid. 
We mention that $\mathcal{N}$ is unoriented or
oriented in the definition of a ported Tutte function
because the big Tutte polynomial never distinguishes
different orientations of the same matroid but
the oriented variant of it, which is defined for
ported oriented matroids, does distinguish some
orientations of the same matroid.

Here is a simple example:  Let the set of ports be $P=\{p_1,p_2\}$.  
Let $\mathcal{N}_1^-(P,\emptyset)$ be the oriented matroid 
with ground set $P$ and oriented
circuit collection $\{++,--\}$. $\mathcal{N}_1^-$ is one of the 
two orientations of the 
rank 1 uniform matroid with ground set $P$.
Let $\mathcal{N}_1^+(P,\emptyset)$ be the other orientation---Its oriented circuit 
collection
is $\{+-,-+\}$.  Since $\mathcal{N}_1^+$ and $\mathcal{N}_1^-$ are 
decomposable under ported Tutte
decomposition, any ported Tutte invariant $F$ of oriented matroids for which
$F(\mathcal{N}_1^+)\neq F(\mathcal{N}_1^-)$ distinguishes different 
orientations of the same matroid.  More interesting examples are given in
\textsection \ref{MaxwellSection}.

The basic theory of non-ported matroid Tutte invariants proves that the
Tutte polynomial has various combinatorial interpretations, i.e., different
expansions over subsets or other structures.  The generalizations of these
for ported Tutte functions and invariants are discussed in \textsection 
\ref{Activity}---They 
all include variables identified with (connected) matroids
on port elements as in Las Vergnas' big Tutte polynomial, except that 
these matroids are oriented.
One such expansion generalizes the corank-nullity polynomial.
One of our results is that our 
extensor-valued Tutte function of ported unimodular oriented matroids
can be expressed by replacing  
each oriented matroid monomial by this function's value, an extensor,
on the oriented matroid which the monomial identifies.
Since the expression also requires $u=v=0$ where $u,v$ are the 
corank-nullity polynomial's variables, we see that our extensor-valued
Tutte function generalizes the basis enumerator.

Our extensor-valued Tutte function provides one 
example of a ported Tutte invariant of graphic oriented matroids (and more
generally, unimodular oriented matroids).  Given a linear representation $N$
of a graphic or other unimodular oriented matroid, the invariant's value
is defined by our extensor function applied to the extensor presenting
$N$, when we take all the parameters $r_e=g_e=1$.  We will see that this
generalizes the fact that the reduced Laplacian determinant equals the number 
of bases in the corresponding graphic matroid, which is the number of 
spanning trees when the graph is connected.

\subsection{Structure of Solutions}

The second involvement between Laplacians, ported Tutte functions and exterior
algebra beyond the matrix tree theorem grows out from the 
discrete (i.e., matrix)
variant of Laplace's equation, the problems in classical physics that it
models, and the structure and solution methods for those problems.
Among many other analogs, the parametrized, discrete equation 
of Laplace combines Kirchhoff's current law (of flow conservation),  Kirchhoff's 
voltage law (electromotive force 
in direct current electrical systems is determined by a potential function),
and Ohm's law (current flow is proportional to potential difference).  
Given a resistive electrical network, the voltages (i.e., potentials)
at all vertices can be determined, by solving Laplace's equation, 
from given voltages at certain vertices and/or the currents into 
other vertices flowing from the environment.
Our work might contribute a few new insights of matroid theory into this 
situation.  One of our starting points for developing this application
is well-known in both matroid theory and electrical network theory:
Certain insights are obtained when one can replace analysis in terms
of graph vertices by analyses involving graph edges and matroids on
them.  Hence, we model the interactions with the environment
by currents through, and voltage drops across port edges, instead of
currents into vertices and potential values at vertices.

Consider all electrical networks with a fixed set of port edges $P$.
Consider two networks to be the same members of this class 
if they have the same edges and the
same cycle spaces.  Such equivalent networks will have the same
electrical behavior.
For us, this class is the same as the class of graphic oriented 
matroids with a fixed set of port edges for which two parameters, $r_e$ and
$g_e$ are given for each non-port edge.  Each network determines
the ported extensor that represents the graph's graphic oriented matroid.
When our extensor function is applied to such an extensor, the
function value is an extensor corresponding to the linear
subspace of solutions projected to the port variables.
For the electrical network application,
the significance of our main theorem 
is that the function of networks with port edges $P$
that gives each network's projected solution space
\textit{is a ported Tutte function.}  Routine network analysis
to solve for port behavior can use the Laplacian determinant and cofactors 
(or equivalent elimination methods) to solve for the constraints between 
port variables.
But our results show that not only is
the Laplacian determinant a Tutte function, but when an extensor
is used to express the solution space, the entire solution
is a Tutte function.
As a corollary, the network solution \emph{when expressed in exterior
algebra} can be written as linear combination of solutions of networks
with \emph{no resistor edges at all---only port edges}.  The coefficients
are homogeneous multilinear products of the parameters.

The extensor value is can be calculated in three ways.  
One is by recursive application of the Tutte equations, i.e., 
Tutte decomposition.  The second is
by substitution into a polynomial that extends the parametrized
Tutte polynomial to ported, parametrized oriented matroids.  
The third is by Gaussian elimination.  
That calculation 
generalizes the evaluation of the Laplacian
determinant.  For graphs it is equivalent to finding a suitable
representation of the solution space of an electrical network
projected onto the variables associated to the port edges.
Of course, our extensor-valued invariant is calculated after setting
all $r_e=g_e=1$.

\subsection{Computational Complexity}

An additional motivation for our work comes from theoretical computer
science.  The number of bases in unimodular matroids is virtually the
only non-trivial Tutte invariant whose computation is tractable.
For such matroids, the computation input is a totally unimodular 
matrix $N$ and the output value is $\det NN^{t}$.
The more general problem, to compute the number of bases a matroid represented
by an arbitrary matrix, is $\#\mathcal{P}$-complete
(see \textsection \ref{Complexity}
for details.)
Our generalization of this $\det NN^{t}$ is an extensor; a succinct
matrix representation of it can be computed by Gaussian elimination (like
a determinant).  Therefore, when the notion of Tutte invariant is
generalized with ports and with exterior algebraic values (and the 
Tutte equations are modified with sign corrections), new computationally 
tractable invariants are obtained.

\subsection{Additional Context}

See \cite{Recski,Murota} for elements of electrical 
network theory from the point of view of matroids.  Electrical networks
with non-linear but monotone resistance functions were studied using
graphs by\cite{HaslerNeirynck,HaslerDApplMath} in a way 
that led us to apply oriented matroids to this topic\cite{sdcOMP}; these
publications may help orient the reader our point of view.
Our current work and \cite{sdcPorted} extend to sets of
more than one distinguished element results and ideas about series and parallel
connection of both networks and matroids.  In particular, the
four-variable Tutte polynomial for a pointed matroid, which was
defined and applied to these topics by Brylawski\cite{BrylawskiPointed}, 
is a special case of our Tutte function.  Our previous
work covered solutions\cite{sdcELEW} by 
energy minimization\cite{LovaszEnergy} of electrical networks
(and generalizations) with multiple ports, applications of
enumeration\cite{sdcISCAS95} combined with oriented matroid
ideas\cite{sdcISCAS98} to network analysis, and a (non-oriented) matroid
abstraction\cite{sdcBDIMatroid} of the solutions.

We proceed to the details about exterior algebra and
\Plucker coordinates pertaining to
realizable non-oriented and oriented matroids.  They include
algebraic operations and identities which correspond to some elementary
matroid relationships.  Our main construction, an extensor valued function
of ported extensors, and the ported Tutte equation variant that it
satisfies, is presented in \textsection \ref{Extensors}.
A variant of the corank-nullity polynomial in \textsection
\ref{Classify} is used to express our function restricted
to unimodular extensors.  The variant differs from
Las Vergnas' big Tutte polynomial 
so (1) it
applies to oriented matroids instead of unoriented matroids, and
(2) it includes parameters as in Definition \ref{firstTutteFunction}.
Extensions to ported matroids of known results
about expressing Tutte functions as set, basis and flat
expansions, and the related open questions follow 
in \textsection \ref{Activity}.
Further discussion of the electrical 
context is 
given in \textsection \ref{ENets}
and brief remarks on peripheral topics appear in
\textsection \ref{Peripheral}.

\section{Preliminaries}
\label{Preliminaries}

Throughout, $\FieldK$ denotes a field, either the 
reals, rationals, or their extensions generated by
the parameters $g_e$, $r_e$.

\subsection{Exterior Algebra}
\label{ExteriorAlgebraSection}

We refer the reader to basic texts such as 
\cite[\textsection 7.1-7.2, on associative and exterior algebras over fields]{JacobsonI}
for complete development and proofs.
The following is a synopsis with the emphasis on the
facts we will need.  It also explains certain notational conventions which
help to mimic oriented matroid theory in exterior algebra.  We use 
a combinatorial approach to adapt the operations of Hodge star
(for duality) and tensor contraction (for matroid contraction).

An associative algebra $\mathcal{A}$ over $\FieldK$ is a ring that is also
a vector space over $\FieldK$, for which addition and $0$ are the
same in both the ring and the vector space, and for which the ring
and scalar multiplications are compatible:  $a(xy)=(ax)y=x(ay)$
where $a\in\FieldK$ and $x,y\in \mathcal{A}$.

Let $V$ be the vector space $\FieldK S$
where finite set $S$ is a basis.  Thus $V$ consists of the
all $\sum_{e\in S}a_e e$, $a_e\in\FieldK$, where
$\sum a_e e = 0$ if and only if $a_e=0$ for every $e\in S$.
The associative algebra over $\FieldK$ generated by $V$ consists
of all finite $\FieldK$-linear combinations of $\mathbf{1}$ (the ring identity) and
formal finite (non-commuting) products of elements of $S$.

The \textbf{exterior algebra} $\mathcal{E}(V)$ over $V$ is the 
quotient of the associative algebra over $\FieldK$ generated by 
$V$ modulo the algebra ideal $I$ generated by products $v^2$, $v\in V$.  
The image of each non-zero $v\in V$ under the map $v\rightarrow v+I\in\mathcal{E}(V)$
is denoted by $\ext{v}$.  These $\ext{v}$ will also be called \textbf{vectors}. 
Thus, for $v_1, v_2\in V$, 
$(\ext{v}_1+\ext{v}_2)(\ext{v}_1+\ext{v}_2)=\ext{0}$,
$\ext{v}_i\ext{v}_i=\ext{0}$ and so $\ext{v}_1\ext{v}_2=-\ext{v}_2\ext{v}_1$ 
in $\mathcal{E}(V)$.  This and
the associativity law imply that each product of a sequence of 
vectors vanishes if the sequence has repeated elements.  Indeed the
product vanishes if and only if there is a linear dependency among
the vectors being multiplied.  Note that a non-zero product of two or more
vectors is not a vector.

A particular basis of $\mathcal{E}(V)$ 
is constructed
from the basis $S=\{s_1, \ldots, s_n\}$ of $V$.  This basis consists
of $\mathbf{1}$ together with the $2^n-1$ products of vectors in 
distinct non-empty subsets of $\{\ext{s}_1, \ldots, \ext{s}_n\}$, each
product written in a particular order.  A formula for exterior
product expressed in terms of this basis is used in \cite{JacobsonI}
to prove that the product is associative.  The formula expresses the
following fact for products of basis vectors which is true for all
products:
Given any sequence of vectors $v_1\cdots v_k$
and permutation $\sigma\in\Perms_k$ with
sign $\epsilon(\sigma)$, the exterior product satisfies the 
\textbf{alternating law}
\[
\ext{v}_1\cdots\ext{v}_k=\epsilon(\sigma)\ext{v}_{\sigma_1}\cdots\ext{v}_{\sigma_k}.
\]
As a result, $\mathcal{E}(V)$ is an associative algebra
that has dimension $2^n$ when viewed as a vector space over $\FieldK$.  
It is customary to use increasing
order of subscripts, so each $X\in\mathcal{E}(V)$ can be expressed by
\[
X = x_{\emptyset}\mathbf{1} + \sum_{\substack{\emptyset\neq A\subseteq S :\\
                           A=\{s_{i_1},\ldots,s_{i_k}\},\;\;  i_1 < \cdots < i_k}}
                         x_A \ext{s}_{i_1}\cdots\ext{s}_{i_k}
\]
with $2^n$ unique coefficients $x_A$, $A\subseteq S$.

We follow a different convention which mimics the one used with
the chirotope cryptomorphism for oriented matroids given in \cite{OMBOOK}.  
For us, $A$ will
denote an arbitrary sequence of elements of $S$.  The
value in $\FieldK$ symbolized by coefficient $x_A$ will depend on the order
as well as the elements of $A$, but these values will satisfy
\[
    x_{A_\sigma}=\epsilon(\sigma)x_A
\]
where $A_{\sigma}=(a_1\cdots a_k)_{\sigma}=a_{\sigma{1}}\cdots a_{\sigma{k}}$
is $A$ permuted by $\sigma$.
In general,
a function like $A\rightarrow x_A$ 
is called \textbf{alternating} if it has this property.
Our convention allows
$A$ to have repeated elements but the alternating property
requires $x_A=0$ for such A.  

We follow a related convention for subset expansions and formulas within them.
When necessary, a set symbol like $A$ in $A\subseteq S$ will denote
distinct elements written in an arbitrary sequence.  But the expansion
or formula will be written only if its value is independent of the sequence
chosen for each symbol.  Furthermore, when $A$ is a sequence of
distinct basis 
vectors, the corresponding product of their images in $\mathcal{E}(V)$
will be denoted by $\mathbf{A}$.  The empty sequence $\emptyset$ corresponds
to $\mathbf{1}\in\mathcal{E}(V)$.  No sequence of linearly independent
vectors corresponds to $\mathbf{0}\in\mathcal{E}(V)$.
The above $2^n$ term basis expansion is thus written 
\[
X = \sum_{\substack{A\subseteq S\\ A=\{a_1,\ldots,a_{\Card{A}}\}}}x_A\ext{a}_1\cdots \ext{a}_{\Card{A}}
  = \sum_{A\subseteq S}x_A\ext{A}.
\]
This expansion follows our convention because
$x_{A_\sigma}\mathbf{A}_\sigma$ $=$ $\epsilon^2(\sigma)x_A\mathbf{A}$ $=$ 
$x_A\mathbf{A}$.  
Note that $X=\mathbf{0}$ if and only if every $x_A=0$.

The concatenation of sequences $A$, $B$, $C$, \ldots is denoted by $ABC\ldots$

\subsection{Extensors with Ground Set}
\label{ExtWithGroundSetSect}

The exterior algebra $\mathcal{E}(V)$ 
is a powerful tool to explore, in a coordinate free way, relationships 
between linear subspaces of $V$\cite{RotaCayley,WhiteCayleyGeoAppl}.  
These relationships are, in other words, the
theorems about the projective geometry whose points are 
the rank 1 (zero-dimensional) subspaces.   The geometric flats 
(the empty set, points, lines, planes, etc., and the whole space) 
correspond to these linear subspaces.  One way to present
a $\FieldK$-realizable matroid is to map each ground set element
to either a point or to the empty flat in this projective geometry.  
The matroid structure is then expressed in terms of incidence of
the images of the ground set elements with the geometric flats.
The formulation is coordinate-free because these
relationships do not change under a change of the basis $S$ for $V$.
We mention this formulation 
to contrast it with our application of exterior algebra.

In our application, each $\FieldK$-realizable matroid with ground
set $S$ will be presented by \emph{a separate} 
(fully decomposable) non-zero
element $\ext{N}\in\mathcal{E}(V)$, where $V=\FieldK S$.  
Each such element $\ext{N}$ 
will determine a linear subspace $L=L(\ext{N})$ of $\FieldK S$.
Consider the family of matrices with 
columns indexed by $S$ whose row space equals $L$.
The matroid
is presented by the linear dependencies among the columns
of any such matrix.  When $L$ is viewed as a linear 
subspace of functions from $S$ to $\FieldK$, the matroid
is the ``function space geometry (or chain-group geometry)
$G(S,L)$'' discussed in \cite[\textsection 1.1.C]{TheoryMatroids}.
As such, each $e\in S$ corresponds to the linear 
functional given by evaluation of $f\in L$ on $e$.
Therefore these functionals, as a
finite subset of the dual space of $L$, 
comprise a vector representation of the matroid.

The members of $L$ present an oriented matroid
$\mathcal{N}(\ext{N})$ by a set of 
\textbf{covectors} $\mathcal{L}$.  Each covector
is the function $l:S\rightarrow\{+,-,0\}$ determined
by some $f\in L$ by $l(e)=\mbox{sign}\big(f(e)\big)$ for
all $e\in S$.  
In other words, if $L$ is viewed as the row space of a matrix,
then each $l\in\mathcal{L}$ is the signature (i.e., the sequence
of $+$, $-$ or $0$s) that indicates the signs in one row $f\in L$.
The signed cocircuits are the 
covectors with minimal non-empty support.
See \cite{BachemKern} for an exposition of
oriented matroids that begins with 
linear subspace presentations including
the cycle and coboundary (or cocycle) spaces 
of graphs.  
Our topic utilizes exterior algebra
and the chirotope given by the
signs of the \Plucker coordinates of $\ext{N}$
to present oriented
matroids.
Theorem \ref{ExtTheory} states the needed details.
Deeper discussions appear in 
\cite[especially \textsection 2.4 on stratifications of the
Grassmann variety and chap. 8 on realizations]{OMBOOK}.

The following
theorem summarizes the facts we will need.  It characterizes those
elements of $\mathcal{E}(V)$ (as fully decomposable) that determine linear subspaces
and present the $\FieldK$-realizable matroids with ground set $S$.
We will call such elements extensors.

\begin{theorem}
\label{ExtTheory}
Given non-zero element $\ext{N}\in\mathcal{E}(V)$ where 
$V=\FieldK S$, the following conditions are equivalent:
\begin{enumerate}
\item
There exist $r$ linearly independent vectors $v_i\in V$ for which
\[
   \ext{N}=\ext{v}_1\cdots \ext{v}_r.
\]
(When $r=0$, the condition is $\ext{N}=\alpha\ext{1}$ for some
$\alpha\in\FieldK$, $\alpha\neq 0$.)

$L$ is the subspace spanned by $\{v_1,\ldots,v_r\}$.

\item
There exists $r$, $0\leq r \leq \Card{S}$, such that the only 
non-zero coefficients $\ext{N}[A]$ in
\[
\ext{N}=\sum_{A\subseteq S}\ext{N}[A]\ext{A}
\]
satisfy $\Card{A}=r$, and the function
$\ext{N}[A]$ from sequences $A$ to $\FieldK$ is alternating and
satisfies the \emph{Grassmann-\Plucker} relationships:

For all length $r$ sequences $A=a_1\cdots a_r$ and $B=b_1\cdots b_r$  over $S$,
\[
\ext{N}[A]\ext{N}\ext[B] = \sum_{i=1}^r 
\ext{N}[b_ia_2\cdots a_r] \ext{N}[b_1\cdots a_i\hat{b_i}\cdots b_r].
\]
(Here, $a_i\hat{b_i}$ means $b_i$ within sequence $B$ is replaced by $a_i$.)

\item There exists a rank $r$ matrix $N$ with $r$ rows and
with columns indexed by
$S$ for which the coefficients of $\ext{N}$ with 
$\Card{A}=r$ satisfy
\[
   \ext{N}[A] = \det N(A)
\]
where $N(A)$ is the submatrix of $N$ with columns $A$, and the other
coefficients are $0$.  (For $r=0$ the condition
is $\ext{N}[\emptyset]\neq 0$ and all other $\ext{N}[A]=0$.)

$L$ is the subspace spanned by the rows of $N$.

If $r\neq 0$ then $N$ and the $v_1,\ldots v_r$ in 
$\ext{N}=\ext{v}_1\cdots\ext{v}_r$
can be chosen so row $i$ of $N$ holds the coefficients for 
writing $v_i$ as a linear combination of vectors from basis $S$.  

\end{enumerate}
\end{theorem}

\begin{definition}[Extensor]
\label{ExtensorDef}
When the conditions in Theorem \ref{ExtTheory} are true about $\ext{N}$, we 
say $\ext{N}$ is decomposable\footnote{We follow \cite{JacobsonI} and other authors who
omit the qualifier ``fully'' in this context.}, and the coefficients denoted by
$\ext{N}[A]$ are called the \Plucker coordinates of $\ext{N}$.
A decomposable element of an exterior algebra is called an 
\textbf{extensor}.
The integer $r$ is called its rank and is denoted by $\rho\ext{N}$. 
\end{definition}

\Remark: We do not define \Plucker coordinates to equivalence classes of
homogeneous coordinates, so 
extensors that differ by a non-zero scalar multiple have different
\Plucker coordinates, even though they represent the 
the same subspace.  

\Remark: $\ext{N}[A]$ is defined as $0$ for all $A$ with 
$\Card{A}\neq\rho\ext{N}$.

The fact that each rank $r$ extensor is the exterior
product of $r$ vectors, and $\ext{v}_1\ext{v}_2=-\ext{v_2}\ext{v}_1$ for 
vectors implies that extensor multiplication satisfies the 
following \textbf{anticommutative law}:
\begin{equation}
\label{Anticommutative}
\ext{N}_1\ext{N}_2=(-1)^{\Rank{\ext{N}_1}\Rank{\ext{N}_2}}
\ext{N}_2\ext{N}_1.
\end{equation}

\begin{theorem}
\label{OrthogSpace}
Suppose $r$, subspace $L$, matrix $N$, extensor $\ext{N}(S)$ and
its \Plucker coordinates are as described in Theorem \ref{ExtTheory}.
\begin{enumerate}
\item There exists a rank $\Card{S}-r$ matrix $N^\perp$ with columns 
indexed by $S$ and with $\Card{S}-r$ rows for which coefficients
with $\Card{A}=r$ satisfy
\[
   \ext{N}[A] = \det N^\perp(\scomp{A})\epsilon(A\scomp{A})
\]
where $\scomp{A}=S\setminus A$ expressed in an arbitrary sequence and
$\epsilon$ is some non-zero alternating sign function of
sequences over $S$.

\item
$L$ consists of all $x\in \FieldK S$ that satisfy the equations
\[
    N^\perp x = 0.
\]
In other words, $L$ (the row space of $N$) and the row space
of $N^\perp$ are orthogonal complements in $\FieldK S$.
\end{enumerate}
\end{theorem}

\begin{proof} 
See \cite[chap. 7]{JacobsonI}.  An elementary proof of 
theorem \ref{OrthogSpace} related to theorem \ref{ExtTheory}
appears in \cite[VII.3 Theorem I]{HodgePedoe1}.  The
equivalence of the Grassmann-\Plucker relationships to
the other conditions is proved in 
\cite[VII.6 Theorem II]{HodgePedoe1}.
\end{proof}

An \textbf{extensor} $\ext{N}$ \textbf{with ground set} $T$ 
is the finite set $T$ paired with an extensor in
$\mathcal{E}(\FieldK T)$.
We use the notation $\ext{N}=\ext{N}(T)$ and $T=S(\ext{N})$ to indicate
that $\ext{N}$ has ground set $T$.  

We need the ground set for the same reason that a ground set is necessary
to define the dual of a matroid with coloops.  Independent sets, or
the collection of sequences $B$ for which chirotope $\chi(B)=+$ is not
sufficient because a loop doesn't appear in any of these objects.  
Furthermore, in our calculations and proofs
we find it very helpful to combine the signs of 
sequences (i.e., permutations of subsets) from several different sets
within one analysis without having to relabel any elements.
Identities like the Tutte equations
relate function values for objects with 
different ground sets.  It is not sufficient for
an extensor or chirotope to be defined up to sign for certain
identities to be valid (not just up to sign).  This validity may
facilitate the use of the identities in computer programs.

\begin{theorem}
Given 
$\ext{N}(S)$, let $N$ be a matrix
satisfying Theorem \ref{ExtTheory}.
\begin{enumerate}
\item
The collection of those $B\subseteq S$ for which $\ext{N}[B]\neq 0$ is the
collection of bases of a matroid with ground set $S$.  

The same matroid is presented by the independent sets of columns of $N$.
\item
The function $\chi$ of sequences over 
$S$ into $\{+1,-1,0\}$ for which
$\chi(B)$ is the sign of $\ext{N}[B]$ is the chirotope function of the
oriented matroid $\mathcal{N}$ denoted by 
$\mathcal{N}=\mathcal{N}(S)=\mathcal{N}(\ext{N})$.

The covectors of $\mathcal{N}$ are presented by the 
signatures of $N$'s row space elements; the signed
circuits (i.e., oriented matroid ``vectors'' with minimal support sets)
are presented by the signatures of the minimal linear dependencies 
among the columns of $N$.

\item
If $\ext{N}[B]\in \{0,\pm 1\}$ 
for all $B$, then $\mathcal{N}(\ext{N})$ is 
the unimodular (or regular) oriented matroid whose chirotope
function satisfies $\chi(B)=\ext{N}[B]$.  Furthermore,
every unimodular oriented matroid can be presented by such an 
$\ext{N}$.

\end{enumerate}
\end{theorem}

\begin{proof}
See \cite{OMBOOK} or \cite[chap. 5]{BachemKern}.  Details pertaining to 
the unimodular matroids including several characterizations are given
in \cite[Theorem 3.1.1, p. 41]{CombinatorialGeometries}.
\end{proof}

\begin{definition}
If $\ext{N}(S)\neq 0$ and $e\in S$ then
\begin{itemize}
\item
$e$ is called a loop if $\ext{N}[B]=0$ for all $B$ with $e\in B$, and
\item
$e$ is called a coloop if every $B$ such that $\ext{N}[B]\neq 0$ 
satisfies $e\in B$.
\end{itemize}
\end{definition}

\Remark: $e$ is therefore a loop or coloop in $\ext{N}$ if and only if
it is a loop or coloop respectively in the matroid presented by 
$\ext{N}$.

\begin{definition}
Each $\ext{N}(S)\neq 0$ defines the function $\rho_{\ext{N}}$ on
subsets $A\subseteq S$ where $\rho_{\ext{N}}(A)$ is the rank of $A$ in 
the matroid presented by $\ext{N}(S)$.
\end{definition}

\begin{theorem}
\label{ExtMinors}
Given $\ext{N}(S)\neq 0$, $e\in S$ and $S'=S\setminus e$:
\begin{enumerate}
\item The \Plucker coordinate function for $\ext{N}$ restricted to sequences
$B\subseteq S'$ is the \Plucker coordinate function for
an extensor denoted by $(\ext{N}\setminus e)(S')$.  
This operation $\ext{N}\rightarrow\ext{N}\setminus e$ is called \textbf{deletion} of $e$.

$(\ext{N}\setminus e)\neq \ext{0}$ if and only if 
$e$ is not a coloop in the matroid presented by $\ext{N}$.
In this case,
the unoriented or oriented matroid minor $\mathcal{N}\setminus e$ 
is presented by $\ext{N}\setminus e$ and
$\rho(\ext{N}\setminus e) = \rho\ext{N}$.  

\item The function defined by $\ext{N}[Be]$ 
for sequences $B\subseteq S'$ is the \Plucker coordinate function for
an extensor denoted by $(\ext{N}/e)(S')$.  
This operation $\ext{N}\rightarrow\ext{N}/ e$ is called 
\textbf{contraction} of  $e$.

$(\ext{N}/e)\neq \ext{0}$ if and only if 
$e$ is not a loop in the matroid presented by $\ext{N}$.
In this case,
the unoriented or oriented matroid minor $\mathcal{N}/ e$ 
is presented by $\ext{N}/e$ and
$\rho(\ext{N}/ e) = \rho\ext{N}-1$.  

\end{enumerate}
\end{theorem}

\Remarks: 
If $\ext{N}=\ext{0}$ then $\ext{N}/e=\ext{N}\setminus e=\ext{0}$.
The zero extensor $\ext{0}$ does not present any matroid.  All rank 0
(empty $S$ or loops only) matroids
have only one basis $\emptyset$; they are presented by the non-zero
extensors $\alpha\ext{1}(S)$ of rank $0$.  

\begin{proof}
Let $N$ be a matrix representing $\ext{N}$ in Theorem \ref{ExtTheory}.

$(\ext{N}\setminus e)(S')$ is the extensor known from Theorem \ref{ExtTheory}
when the column labeled by $e$ is deleted from $N$.  Note that 
if $e$ is a coloop then
this reduces the rank of $N$ and so $(\ext{N}\setminus e)(S')=\ext{0}$.  

If $e$ is a loop in $\ext{N}$ then $(\ext{N}/e)=\ext{0}$.  Otherwise, 
the $\ext{N}[Be]$ are fixed non-zero multiples 
of the \Plucker coordinates from a 
matrix obtained from the $N$ by row operations to make all
but one entry in column $e$ zero and then deleting the row and column
with that non-zero entry.

See \cite[\textsection 3.5]{OMBOOK} for oriented matroid minors and other
structures in terms of chirotopes.
\end{proof}

\begin{theorem}
\label{ExtDC}
Given $\ext{N}(S)$ and $e\in S$,
\[
\ext{N}(S) = (\ext{N}/e)\ext{e} + (\ext{N}\setminus e)\ext{1}(\{e\})
\]
\end{theorem}
The multiplication by $\ext{1}(\{e\})$ makes the ground set of the 
second term be $S$ instead of $S\setminus e$.  It will be omitted in
contexts where the ground set is clear.

\begin{proof}
Let $B\subseteq S$.
We prove that each \Plucker coordinate $\ext{N}[B]$ equals the sum of the 
corresponding \Plucker coordinates of the extensors on the right.

If $e\in B$ we can write $B=B'e$.
$\ext{N}[B'e]$ $=$ $(\ext{N}/e)[B']$ 
$=$ $(\ext{N}/e)\ext{e}[B'e]$, and $(\ext{N}\setminus e)[B]=0$.

If $e\not\in B$ then $(\ext{N}/e)\ext{e}[B]=0$
and $\ext{N}[B]=(\ext{N}\setminus e)[B]$.
\end{proof}

It is convenient to let $\ext{N}/A$ denote 
$\ext{N}/a_k/\cdots/a_1$ where $A=a_1\cdots a_k$, and similarly for
$\ext{N}\setminus A$.  It follows that $\ext{N}/A[X]$
$=$ $\ext{N}[XA]$ for all $X$.  We note that for $\sigma\in\Perms_k$,
\begin{equation}
\label{PermuteContraction}
\begin{split}
\ext{N}/A_\sigma &= \epsilon(\sigma)\ext{N}/A\text{, but}\\
\ext{N}\setminus A_\sigma &= \ext{N}\setminus A.
\end{split}
\end{equation}

\subsection{Ground Set Orientation and Duality}

\begin{definition}[Ground set orientation]
An orientation of the ground set $\epsilon$ is an alternating function
into $\{+1, -1, 0\}$ of sequences of ground set elements that is 
non-zero on sequences of distinct elements, and which satisfies
$\epsilon(\emptyset)=1$.
\end{definition}

One family of ground set orientations is derived from fixed linear
orders on all ground set elements using the rule that
$\epsilon(X)=(-1)^v$ where $v$ is the number of inversions in $X$
(where an inversion is $(i,j)$ with $i<j$ and $x_i>x_j$).  A permutation 
$\sigma\in\Perms_n$ of $\{1,\ldots,n\}$ is always considered a sequence
$\sigma_1\sigma_2\ldots\sigma_n$ of natural numbers 
with ground set orientation derived from their usual ordering.
Hence, $\epsilon(\sigma)$ is the usual sign of permutation $\sigma$.
However, ground set orientations of matroid elements or graph edges
will \emph{not} be assumed to derive from a linear order. 

Since permutations $\sigma\in\Perms_n$ and sequences of ground set
elements will not be confused, we use the same symbol $\epsilon$
for permutation sign and ground set orientation.

Given a length $n$ sequence $X=x_1\ldots x_n$ and $\sigma\in\Perms_n$,
let $X_\sigma$ denote $x_{\sigma_1}\ldots x_{\sigma_n}$.
The following routine facts will be used in our proofs:  Of course,
$F$ is alternating means $F(X_\sigma)$ $=$ $\epsilon(\sigma)$
$F(X)$ for all sequences $X$ and $\sigma\in\Perms_{\Card{X}}$.
\begin{lemma}
\label{lemmaEpsilonFacts}
Suppose $\epsilon_1$ and $\epsilon_2$ are arbitrary alternating 
functions of sequences.
\begin{enumerate}
\item If $n=|X|=|Y|$, $\sigma$ $\in \Perms_n$,
and $A, X, B$, $C, Y, D$ are sequences then
\begin{eqnarray*}
\epsilon_1(AXB)\epsilon_2(CYD)=&
\epsilon(\sigma)\epsilon_1(AX_\sigma B)\epsilon_2(CYD) \\
                              =&
\epsilon(\sigma)\epsilon_1(AXB)\epsilon_2(CY_\sigma D) \\
                              =&
\epsilon_1(AX_\sigma B)\epsilon_2(CY_\sigma D)
\end{eqnarray*}

\item 
$\epsilon_i(XY) =
(-1)^{|X||Y|}\epsilon_i(YX)$.
\end{enumerate}
\end{lemma}

With a ground set orientation $\epsilon$ in hand, we define: 
\begin{definition}[Canonical Dual]
\label{definitionCanonicalDual}
Given $\ext{N}(S)$, 
$\ext{N}^\perp[X] = \ext{N}^{\perp_{\epsilon}}[X]=\ext{N}[\scomp{X}]\epsilon(\scomp{X}\scma X)$, where $\scomp{X}$ is any sequence of the distinct elements
in $S\setminus X$.

The symbol $\perp_{\epsilon}$ will be abbreviated by $\perp$ when 
$\epsilon$ is irrelevant or doesn't require emphasis.
\end{definition}

\Remark:
Each ground set $S$ determines with $\epsilon$ 
a sign choice from among the two 
that both provide a presentation of the
oriented matroid dual.

The demonstration in \cite[end of \textsection 3.6]{OMBOOK}
of oriented matroid chirotope dualization 
has a similar formula, whose right hand side is
independent of an arbitrarily chosen sequence.  It follows that 
our extensor dualization corresponds to the oriented
matroid dualization of the oriented matroid presented by
$\ext{N}$.  Theorem \ref{OrthogSpace} justifies this for
realizable matroids.  (Dualization is also the Hodge star 
operator\cite{HodgePedoe1}
when $S$ is identified with the corresponding basis of the dual 
space. Also see \textsection\ref{GroundSetOrientationBack}.)

\subsection{Identities}
\label{IdentitiesSect}

Our main proof uses some identities on extensors 
that correspond to well-known relationships among matroid operations.
These identities involve extensors with ground sets
for which a ground set orientation is used to define dualization.
The union of disjoint sets is denoted by $\dunion$.

\begin{theorem}
\label{Identities}
\item
\begin{equation}
\label{LinearityPerp}
\begin{split}
(\ext{N}_1+\ext{N}_2)^\perp = \ext{N}_1^\perp + \ext{N}_2^\perp\\
(\alpha\ext{N})^\perp = \alpha\ext{N}^\perp
\end{split}
\end{equation}

\item
\begin{equation}
\label{PerpPerp}
\ext{N}^{\perp\perp}(S) = (-1)^{\rho\ext{N}\;(\Card{S}-\rho\ext{N})}\;\;\ext{N}(S)
                        = (-1)^{\rho\ext{N}\;\rho\ext{N}^\perp}\;\;\ext{N}(S)
\end{equation}

\item
Given $\ext{N}(S)$, and sequences $X\subseteq S$ and $S'=S\setminus X$,
\begin{equation}
\label{DeletePerp}
(\ext{N}\setminus X)^\perp = \epsilon(S')\epsilon(S'X)\;\;(\ext{N}^\perp/X)
\end{equation}
\begin{equation}
\label{ContractPerp}
(\ext{N}/X )^\perp = \epsilon(S')\epsilon(S'X)
(-1)^{\Card{X}\;(\Card{S}-\rho\ext{N})}\;\;(\ext{N}^\perp\setminus X)
\end{equation}

\item 
Given $\ext{N}_i(S_i)$ with $S_1\cap S_2=\emptyset$, the extensor
product $\ext{N}_1\ext{N}_2(S_1\dunion S_2)$ presents the (oriented)
matroid direct sum and 
\begin{equation}
\label{ProdPerp}
(\ext{N}_1\ext{N}_2)^\perp 
= \epsilon(S_1)\epsilon(S_2)\epsilon(S_1S_2)
(-1)^{\rho\ext{N}_1^\perp\rho\ext{N}_2}\;\ext{N}_1^\perp\ext{N}_2^\perp
\end{equation}
\end{theorem}

\begin{proof}
Linearity (\ref{LinearityPerp}) is immediate from the \Plucker coordinate
definition.  It will be repeatedly used with $\alpha = \pm 1$ below.

To prove Theorem \ref{Identities} (\ref{PerpPerp}), write
\begin{equation*}
\ext{N}^{\perp\perp}[A] = (\ext{N}^\perp)^\perp[A] 
                  = (\ext{N}^\perp)[\scomp{A}]\;\epsilon(\scomp{A} A),\\
= \ext{N}[\scomp{\scomp{A}}]\;
\epsilon(\scomp{\scomp{A}}\;\scomp{A})\epsilon(\scomp{A} A)\\
= \ext{N}[A]\;
\epsilon(A\scomp{A})\epsilon(\scomp{A} A)
\end{equation*}
where in the last equation we chose the sequence order $\scomp{\scomp{A}}=A$.
Therefore the sign correction is 
$(-1)^{|A|\;|\scomp{A}|}$.
For 
non-zero coordinates this is 
$(-1)^{\rho\ext{N}\;(\Card{S}-\rho\ext{N})}$ $=$
$(-1)^{\rho\ext{N}\;\rho\ext{N}^\perp}$.

To prove Theorem \ref{Identities} (\ref{DeletePerp}), write
\begin{equation*}
(\ext{N}\setminus X)^\perp[A] = 
(\ext{N}\setminus X)[\scomp{A}]\;\epsilon(\scomp{A}A)=
                          \ext{N}[\scomp{A}]\;\epsilon(\scomp{A}A),
\end{equation*}
where $\scomp{A}=S'\setminus A$.  But in 
\begin{equation*}
(\ext{N}^\perp/X)[A] = \ext{N}^\perp[AX] = \ext{N}[\scomp{AX}]\;
                  \epsilon(\scomp{AX}AX)
\end{equation*}
the elements in sequence 
$\scomp{AX}$ 
are $S'\setminus A$, the
same as in the sequence symbolized by $\scomp{A}$ in the previous equation.
We can therefore choose $\scomp{AX}$ $=$ $\scomp{A}$ and write
\begin{equation*}
(\ext{N}^\perp/X)[A] = \ext{N}[\scomp{A}]\;
                  \epsilon(\scomp{A}AX).
\end{equation*}
Combining the two sign corrections gives 
$\epsilon(\scomp{A}A)\epsilon(\scomp{A}AX)$.  That equals
$\epsilon(S')\epsilon(S'X)$ for all reorderings
$S'$ of $\scomp{A}A$.

We can get (\ref{ContractPerp}) from (\ref{PerpPerp}) and (\ref{DeletePerp}).
Specifically, $(\ext{N}/X)^\perp$ $=$ 
$(\ext{L}^\perp/X)^\perp(-1)^{\Rank{\ext{N}}\Rank{\ext{N}^\perp}}$ with
$\ext{L}=\ext{N}^\perp$.  This equals
\[
(\ext{L}\setminus X)^{\perp\perp}\;
   \epsilon(S')\epsilon(S'X)(-1)^{\Rank{\ext{N}}\Rank{\ext{N}^\perp}}
\]
\[
= (\ext{N}^\perp\setminus X)\;
(-1)^{\rho(\ext{N}^\perp\setminus X)\rho(\ext{N}/X)}
   \epsilon(S')\epsilon(S'X)(-1)^{\Rank{\ext{N}}\Rank{\ext{N}^\perp}}
\]
As usual, we can restrict attention to non-zero coordinates.
Let $\rho \ext{N}=r$, $\Card{S}=s$, and $\Card{X}=x$ so $\rho(\ext{N}/X)=r-x$
and $\rho(\ext{N}^\perp\setminus X)$ $=$ $\rho\ext{N}^\perp$ $=$ $s-r$. The sign
correction is therefore
\[
   \epsilon(S')\epsilon(S'X)(-1)^{(s-r)(r-x)+r(s-r)}=
   \epsilon(S')\epsilon(S'X)(-1)^{(s-r)(2r-x)}=   
   \epsilon(S')\epsilon(S'X)(-1)^{(s-r)x}.
\]

To prove Theorem \ref{Identities} (\ref{ProdPerp}), take
$\scomp{A_1A_2}=\scomp{A_1}\;\scomp{A_2}$ with each 
$\scomp{A_i}=S_i\setminus A_i$ in
\begin{equation*}
\begin{split}
(\ext{N}_1 \ext{N}_2)^\perp[A_1 A_2]=\\
(\ext{N}_1 \ext{N}_2)[\scomp{A_1} \scomp{A_2}]\;
  \epsilon(\scomp{A_1} \scomp{A_2} A_1 A_2 )\\
=\ext{N}_1[\scomp{A_1}]\; \ext{N}_2 [\scomp{A_2}]\;
  \epsilon(\scomp{A_1}\;\scomp{A_2} A_1 A_2 )\\
=\ext{N}_1^\perp[A_1]\; \epsilon(\scomp{A_1}A_1)\;
 \ext{N}_2^\perp[A_2]\; \epsilon(\scomp{A_2}A_2)\;
  \epsilon(\scomp{A_1} (\scomp{A_2} A_1) A_2 )\\
=\ext{N}_1^\perp[A_1]\;  \ext{N}_2^\perp[A_2]\;
  \epsilon(\scomp{A_1}A_1)\epsilon(\scomp{A_2} A_2)
  \epsilon(\scomp{A_1} (A_1 \scomp{A_2})  A_2 )
  (-1)^{\Card{A_1}\Card{\scomp{A_2}}}\\
=\ext{N}_1^\perp[A_1]\;  \ext{N}_2^\perp[A_2]\;
  \epsilon(S_1)\epsilon(S_2)
  \epsilon(S_1 S_2)
  (-1)^{\rho\ext{N}_1^\perp\rho\ext{N}_2}
\end{split}
\end{equation*}
where in the last equation, we applied permutations $\sigma$
and $\tau$, each twice, for which $(\scomp{A_1}A_1)_\sigma=S_1$
and $(\scomp{A_2}A_2)_\tau=S_2$.  We then substituted the 
correct ranks for cases where the coordinate is not 0.
\end{proof}

\section{An  Extensor Tutte Function}
\label{Extensors}

Recall that a ported extensor or matroid is one whose ground set
has a distinguished subset of port elements.

Given a ported extensor $\ext{N}(P,E)$ (the notation means $P$ is
the set of ports and the ground set is $P\dunion E$), 
we will define a parametrized extensor $\ext{M}_E(\ext{N})$ 
using extensor operations.  We will illustrate its construction
with extensors and equivalent matrices.
We then write and prove
parametrized identities
satisfied by the function $\ext{N}(P,E)\rightarrow \ext{M}_E(\ext{N})$ which 
are analogous to the
ported Tutte equations.  
Our identities however apply to extensors rather than to commutative ring
values.  The identities include
\textit{sign-correction factors} that depend on the particular ground set
orientation $\epsilon$ that was used to define 
$\ext{M}_E(\ext{N})$.  

The definition of $\ext{M}_E(\ext{N})$ below applies to all
extensors $\ext{N}$ over $\FieldK (E\dunion P)$.  The main
result therefore belongs to exterior algebra.
Section \ref{Classify} shows how $\ext{M}_E$ defines the extensor valued 
function on the minor closed class of ported unimodular oriented matroids
$\mathcal{N}(P,E)$ by 
$\ext{M}_E(\mathcal{N})=\ext{M}_E(\ext{\pm N})$ where $\mathcal{N}$
is presented by either unimodular extensor $\pm \ext{N}$.  

Each linear map on $V$ can be extended to a unique exterior algebra 
map on the exterior algebra $\mathcal{E}(V)$
\cite[Theorem 7.1]{JacobsonI}.  Given ground set $P\dunion E$, let
$P_\Vs$ and $P_\Is$ be two disjoint copies of $P$, both also
disjoint from $E$.  For each $p\in P$,
let $p_\Vs\in P_\Vs$ and $p_\Is\in P_\Is$ be the corresponding elements 
in the respective copies.  We define the following maps
from $\FieldK (P\dunion E)$ to $\FieldK[r_e,g_e] (P_\Vs\dunion P_\Is\dunion E)$
and extend them to the exterior algebra.  If necessary,
the field $\FieldK$ is extended with the parameters $g_e,r_e$,
$e\in E$.
\begin{equation}
\begin{split}
\Vs_r(\ext{e}) = r_e \ext{e} \text{ for $e\in E$ and }
\Vs_r(\ext{p}) = \ext{p}_\Vs \text{ for $p\in P$.}\\
\Is_g(\ext{e}) = g_e \ext{e} \text{ for $e\in E$ and }
\Is_g(\ext{p}) = \ext{p}_\Is \text{ for $p\in P$.}
\end{split}
\end{equation}

In terms of matrices, $\Is_g$ signifies multiplying
column labeled $e$ by $g_e$ for each $e\in E$ and renaming
column $p$ by $p_\Is$ for each $p\in P$.  Likewise, $\Vs_r$
signifies multiplying
column $e$ by $r_e$ and renaming column $p$ by $p_\Vs$.

The parameter subscript symbols in $\Is_g$ and $\Vs_r$ will sometimes be
omitted for brevity.  For subset $Q\subseteq P$, 
$Q_{\Vs}$ denotes $\{q_{\Vs} : q\in Q\}\subseteq P_{\Vs}$ and 
$Q_\Is$ denotes $\{q_{\Is} : q\in Q\}\subseteq P_{\Is}$.
Recall that set symbols denote sequences.  The sequences of
$Q_{\Vs}$ and $Q_\Is$ correspond to the sequence of $Q$.

\begin{definition}
\label{DefM}
Given 
a ported extensor \(\ext{N}(P,E)\),
a ground set orientation \(\epsilonfun\) and dual operator
$\perp_\epsilon$,
parameters \(g_e\) and \(r_e\) for each 
\(e\in E\), and 
\(\epsilonfun\)-preserving functions $\Vs_r$ and $\Is_g$ defined 
above, let
\begin{equation}
\begin{split}
\label{MENDef}
\ext{M}(\ext{N}) = \Is_{g}(\ext{N})\;\Vs_{r}(\ext{N}^{\perp_\epsilon})
\text{ and }\\
\ext{M}_E(\ext{N}) =
 \ext{M}(\ext{N})/E
\end{split}
\end{equation}
\end{definition}

Hence, $\ext{M}_E(\ext{N})$ is defined as a ported extensor
$\ext{M}(\ext{N})=\ext{M}(\ext{N})(P_\Is\dunion P_\Vs,E)$ 
contracted by the sequence 
of non-port elements $E$.  Therefore, 
by Theorem \ref{ExtMinors}
it is an extensor.  
Each pair of sequences
$I\subseteq P$, $V\subseteq P$ with $|I|+|V|=|P|$
specifies the \Plucker coordinate
of $\ext{M}_E(\ext{N})$ 
with index $I_\Is V_\Vs$ and value 
\begin{equation}
\label{MENComponentDef}
\ext{M}_E(\ext{N})[I_{\Is}V_{\Vs}] 
      = \big(\Is(\ext{N})\;\;\Vs(\ext{N}^\perp)\big)[I_{\Is}V_{\Vs}E].
\end{equation}

\begin{proposition} 
\label{PropostionAlphaSquared}
For $\alpha\in\FieldK$, 
$\ext{M}_E(\alpha \ext{N}) = \alpha^2 \ext{M}_E(\ext{N})$.
\end{proposition}
\begin{proof} 
$\ext{M}(\alpha \ext{N}) = \alpha^2 \ext{M}(\ext{N})$ is immediate
from the definition. Contraction $\ext{M}/E$ is linear in 
$\ext{M}$.
\end{proof}

We can express  (\ref{MENDef}) in matrix terms.
Let $N$ be some full row rank matrix with columns indexed
by $P\dunion E$ that presents $\ext{N}(P\dunion E)$.
Similarly, let $N^{\perp}$ denote a matrix presentation
of $\ext{N}^{\perp}$.

\begin{figure}
\begin{center}\input{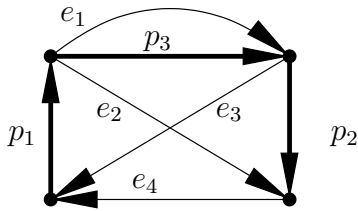}\end{center}
\caption{Graph defining the graphic oriented matroid $\mathcal{N}$}
\label{K4figure}
\end{figure}

\textbf{Example.} 
We show one totally unimodular matrix
representation $N$ of the 
ported graphic oriented matroid with $P=\{p_1,p_2,p_3\}$
and $E=\{e_1,e_2,e_3,e_4\}$ for the graph in figure \ref{K4figure}.
The rows code 3 oriented cutsets 
which determine 
a basis for the 1-coboundary (or cocycle) space.  We also
express $\ext{N}$ by the exterior 
product 
of the vectors given by the rows of this matrix.


\[
N=\begin{array}{c}
\begin{array}{ccccccc}
p_1 & p_2 & p_3 & e_1 & e_2 & e_3 & e_4 \\[-1em]
\hphantom{+1} & \hphantom{+1} &\hphantom{+1} &\hphantom{+1} &\hphantom{+1} &\hphantom{+1} &\hphantom{+1} 
\end{array} \\ 
\left[\begin{array}{ccc|cccc}
 -1 &  0  &  +1 & +1  &  +1 &  0  & 0   \\ 
  0 & +1  &  -1 & -1  &   0 & +1  & 0   \\ 
 -1 & -1  &  +1 & +1  &   0 &  0  &+1  
\end{array}\right]
  \end{array}
\]

\[
\ext{N}=
\begin{array}{c}
(-\ext{p}_1+\ext{p}_3+\ext{e}_1+\ext{e}_2)\cdot\\
(\ext{p}_2-\ext{p}_3-\ext{e}_1+\ext{e}_3)\cdot\\
(-\ext{p}_1-\ext{p}_2+\ext{p}_3+\ext{e}_1+\ext{e}_4)
\end{array}
\]

Next, we write one totally unimodular 
matrix $N^\perp$ for the canonical
dual.  We have checked that the sign satisfies Definition 
\ref{definitionCanonicalDual} with
$\epsilonfun\ $ chosen so $\epsilon(p_1p_2p_3e_1e_2e_3e_4)=1$
by verifying 
$\ext{N}^\perp[e_1e_2e_3e_4]$ 
$=$
$\ext{N}[p_1p_2p_3]\epsilon(p_1p_2p_3e_1e_2e_3e_4)$.


\[
N^\perp=
\begin{array}{c}
\begin{array}{ccccccc}
p_1 & p_2 & p_3 & e_1 & e_2 & e_3 & e_4 \\[-1em]
\hphantom{+1} & \hphantom{+1} &\hphantom{+1} &\hphantom{+1} &\hphantom{+1} &\hphantom{+1} &\hphantom{+1} 
\end{array} \\ 
\left[\begin{array}{ccc|cccc}
  0 &  0  &  +1 &  -1 &   0 &  0  & 0 \\
 +1 &  +1 &  +1 &   0 &   0 &  0  & +1\\
  0 &  +1 &  +1 &   0 &  -1 &  0  &  0\\
 +1 &   0 &  +1 &   0 &   0 & +1  &  0
	\end{array}\right]
  \end{array}
\]

Continuing the general discussion, let 
$G$ and $R$ be the diagonal matrices of the $g_e$ and $r_e$.
The matrix
\begin{equation}
\label{matrixM}
M(N)
=
\left[
	\begin{array}{ccc}
		N(P) & 0 & N(E)G \\
		0    & N^{\perp}(P) & N^{\perp}(E)R
	\end{array}
  \right]
\end{equation}
has order $(p+e)\times(2p+e)$, columns indexed by sequence
$P_{\Is}P_{\Vs}E$ it and presents $\ext{M}(\ext{N})$.

\textbf{Example continued.}
We abbreviate labels $p_{\Is 1}$ and $p_{\Vs 1}$ by $i_1$ and $v_1$, etc.


\[
M(N)=
\begin{array}{c}
\begin{array}{cccccccccc}
i_1 & i_2 & i_3 & v_1 & v_2 & v_3  & e_1  & e_2 & e_3 & e_4 \\[-1em]
\hphantom{+1}&\hphantom{+1}&\hphantom{+1}&
\hphantom{+1}&\hphantom{+1}&\hphantom{+1}&
\hphantom{+g_1}&\hphantom{-g_1}&\hphantom{g_3}&\hphantom{g_4}
\end{array}\\
\left[\begin{array}{ccc|ccc|cccc}
 -1 &  0  &  +1 &   0 &  0  &   0  & g_1  & g_2 &  0  & 0   \\ 
  0 & +1  &  -1 &   0 &  0  &   0  & -g_1 &   0 & g_3 & 0   \\ 
 -1 & -1  &  +1 &   0 &  0  &   0  & g_1  &   0 &  0  & g_4   \\ \hline
  0 &  0  &   0 &  0  &  0  &  +1  &-r_1  &   0 &  0  & 0 \\
  0 &  0  &   0 & +1  &  +1 &  +1  &   0  &   0 &  0  & r_4\\
  0 &  0  &   0 &  0  &  +1 &  +1  &   0  &-r_2 &  0  &  0\\
  0 &  0  &   0 & +1  &   0 &  +1  &   0  &   0 & r_3 &  0
\end{array}\right]
\end{array}
\]

Generally, $\ext{M}(\ext{N})(P_{\Is}\dunion P_{\Vs}\dunion E)$
is the exterior product of the vectors in
$\FieldK (P_{\Is}\dunion P_{\Vs}\dunion E)$ corresponding to the
rows of this matrix.
$\ext{M}_E(\ext{N})(P_{\Is}\dunion P_{\Vs})$ appears
in the expression
\[
\ext{M}(\ext{N})=\big(\ext{M}_E(\ext{N})\big)\ext{E} + \cdots
\]
where the initial term is the only one with factor $\ext{E}$.

\textbf{Example continued.}
We calculate $\ext{M}_E(\ext{N})$ by doing ring operations
on rows to eliminate all but one non-zero entry in each $E$ column in $M(N)$.  
The result is that
\[
g_1 g_2 g_3 g_4 r_1^6 r_2 r_3 r_4  \ext{M}(\ext{N})
\]
is equal to the following extensor in matrix form:
\[\hspace{-7pt}
\begin{array}{c}
\begin{array}{cccccccccc}
i_1 & i_2 & i_3 & v_1 & v_2 & v_3  & e_1  & e_2 & e_3 & e_4 \\[-1em]
\hphantom{-r_1r_4}&\hphantom{-r_1r_4}&\hphantom{-r_1r_4}&
\hphantom{-g_4r_1}&\hphantom{-g_4r_1}&\hphantom{-g_1r_3-g_3r_1}&
\hphantom{-g_1r_1}&\hphantom{-g_2r_1r_2}&\hphantom{g_3r_1r_3}&\hphantom{g_4r_1r_4}
\end{array}\\
\left[\begin{array}{ccc|ccc|cccc}
-r_1r_2 & 0      &r_1r_2  &  0    &g_2r_1 &g_1r_2+g_2r_1  &  0    &  0       &  0      &  0      \\
0       &r_1r_3  &-r_1r_3 &-g_3r_1&  0    &-g_1r_3-g_3r_1&  0    &  0       &  0      &  0      \\
-r_1r_4 &-r_1r_4 &r_1r_4  &-g_4r_1&-g_4r_1&g_1r_4-g_4r_1 &  0    &  0       &  0     &  0      \\ \hline
0       &0       &0       &0      &0      & g_1          &-g_1r_1&  0       &  0      &  0      \\ 
0       &0       &0       &g_4r_1 &g_4r_1 & g_4r_1       &  0    &  0       &  0      &g_4r_1r_4\\ 
0       &0       &0       &0      &g_2r_1 & g_2r_1       &  0    &-g_2r_1r_2&  0      &  0      \\ 
0       &0       &0       &g_3r_1 &0      & g_3r_1       &  0    &  0       &g_3r_1r_3&  0      
\end{array}\right]\end{array}
\]
After some cancellation, we can read off the answer from the $3\times 6$ upper left submatrix, which is a matrix presentation of the
extensor $r_1^2\ext{M}_E(\ext{N})$:
\[
\begin{array}{c}
\begin{array}{cccccc}
i_1 & i_2 & i_3 & v_1 & v_2 & v_3 \\[-1em]
\hphantom{-r_1r_4}&\hphantom{-r_1r_4}&\hphantom{-r_1r_4}&
\hphantom{-g_4r_1}&\hphantom{-g_4r_1}&\hphantom{-g_1r_3-g_3r_1}
\end{array}\\
\left[
\begin{array}{ccc|ccc}
-r_1r_2 & 0      &r_1r_2  &  0    &g_2r_1 &g_1r_2+g_2r_1     \\
0       &r_1r_3  &-r_1r_3 &-g_3r_1&  0    &-g_1r_3-g_3r_1   \\
-r_1r_4 &-r_1r_4 &r_1r_4  &-g_4r_1&-g_4r_1&g_1r_4-g_4r_1  
\end{array}\right]\end{array}
\]

\Remark:  Each
\Plucker coordinate of $\ext{M}_E(\ext{N})$ is a homogeneous polynomial
of degree
$\Card{E}$ in the $g_e$, $r_e$.  However, this example demonstrates
that there sometimes doesn't exist a matrix expression for 
$\ext{M}_E(\ext{N})$ all of whose entries are polynomials.  The reader can
verify that each order 3 minor of the above matrix is a multiple of
$r_1^2$.

\textbf{Example continued:}
Here are examples of \Plucker coordinates, which can be calculated from
the above matrix as order 3 minors divided by $r_1^2$.
See \textsection \ref{MaxwellSection}.
\begin{eqnarray*}
\ext{M}_E(\ext{N})[v_1v_2v_3]=&g_1g_2g_3r_4+g_1g_2g_4r_3+g_1g_3g_4r_2+g_2g_3g_4r_1 \\
\ext{M}_E(\ext{N})[i_1v_2v_3]=&(g_1r_3+g_3r_1)(g_2r_4+g_4r_2) \\
\ext{M}_E(\ext{N})[v_1i_1v_3]=&-g_1g_4r_2r_3+g_2g_3r_1r_4
\end{eqnarray*}

Observe $\ext{M}_E(\ext{N})[v_1v_2v_3]$ 
is the basis enumerator for $\mathcal{N}(\ext{N})\setminus P$.
The graph and extensor minors corresponding to the terms
of $\ext{M}_E(\ext{N})[v_1i_1v_3]$ are shown in figure \ref{2MinFig}.
\begin{figure}
\[
\begin{array}{ccc}
\text{Graph Minor} & \text{Extensor Minor} & 
\text{Term in } \ext{M}_E(\ext{N})[v_1i_1v_3] \\[.5em]
\parbox[c]{1in}{\input{opp1.pstex_t}} & \ext{N}/\{e_1,e_4\}|P & -g_1g_4r_2r_3    \\
\parbox[c]{1in}{\input{opp2.pstex_t}} & \ext{N}/\{e_2,e_3\}|P & +g_2g_3r_1r_4 
\end{array}
\]
\caption{The two graph and extensor minors with corresponding terms
in $\ext{M}_E(\ext{N})[v_1i_1v_3]$.}
\label{2MinFig}
\end{figure}

\subsection{Main Result}

\begin{theorem}
\label{TheoremMain}
The parametrized extensor-valued function 
$\ext{M}_E(\ext{N})(P_\Vs\dunion P_\Is)$ of ported extensor
$\ext{N}=\ext{N}(P,E)$
has the following properties:
\begin{enumerate}
\item
Given $\ext{N}_1(P_1,E_1)$ 
and $\ext{N}_2(P_2,E_2)$ with $E=E_1\dunion E_2$ and 
$P=P_1\dunion P_2$,
\begin{equation}
   \ext{M}_E(\ext{N}_1\extvee\ext{N}_2)(P,E) =
   \epsilon(P_1P_2E)\epsilon(P_1E_1)\epsilon(P_2E_2)\;
   \ext{M}_{E_1}(\ext{N}_1)\extvee
   \ext{M}_{E_2}(\ext{N}_2).
\end{equation}

\item
If $e\in E$ and $E'=E\setminus e$ then
\begin{equation}
   \ext{M}_E(\ext{N}) = \\
   \epsilon(PE)\epsilon(PE')\;
   \big(
   g_e\ext{M}_{E'}(\ext{N}/ e)
   +
   r_e\ext{M}_{E'}(\ext{N}\setminus e)\big).
\end{equation}
\item
Let $E=\emptyset$.  The \Plucker coordinates of 
$\ext{M}_\emptyset(\ext{N})(P_\Is\dunion P_\Vs)$ satisfy
\[
   \ext{M}_\emptyset(\ext{N})[I_{\Is}V_{\Vs}]=\ext{M}[I_{\Is}V_{\Vs}]=
    \epsilon(\scomp{V}\scma V)\;\ext{N}[I]\ext{N}[\scomp{V}].
\]
for all $I\subseteq P$ and $V\subseteq P$.

\item $\ext{M}_E(\ext{0}) = \ext{0}$.

\end{enumerate}
\end{theorem}

Properties 1.\  and 2.\  in terms of \Plucker coordinates are:
\begin{eqnarray}
\label{Thm10Pt1LHS}
\nonumber   \ext{M}_E(\ext{N}_1\extvee\ext{N}_2)[I_{1I}V_{1V}I_{2I}V_{2V}] =\\ 
   \epsilon(P_1P_2E)\epsilon(P_1E_1)\epsilon(P_2E_2)\;
   \ext{M}_{E_1}(\ext{N}_1)[I_{1I}V_{1V}]\;
     \ext{M}_{E_2}(\ext{N}_2)[I_{2I}V_{2V}].
\end{eqnarray}
\begin{eqnarray}
\nonumber   \ext{M}_E(\ext{N})[I_\Is V_\Vs]  = \\
   \epsilon(PE)\epsilon(PE')\;
   \big(
   (g_e\ext{M}_{E'}(\ext{N}/ e)[I_\Is V_\Vs]
   +
   r_e\ext{M}_{E'}(\ext{N}\setminus e)[I_\Is V_\Vs]\big).
\end{eqnarray}

\Remarks:  
\begin{enumerate}
\item 
Property 2. implies that every
linear combination $g_e\ext{M}_{E'}(\ext{N}/e)$ $+$
$r_e\ext{M}_{E'}(\ext{N}\setminus e)$ is decomposable, i.e., an 
extensor.
\item Proposition \ref{PropostionAlphaSquared} with $\alpha=\pm 1$
implies $\ext{M}_E(\ext{N}_1\ext{N}_2)=\ext{M}_E(\ext{N}_2\ext{N}_1)$.  
We can also verify this from the right hand side of 
property 1 using Lemma (\ref{lemmaEpsilonFacts}.3) and noting the rank of
extensor $\ext{M}_{E_i}(\ext{N}_i)$ is $\Card{P_i}$.
\item
If $\ext{N}\neq \ext{0}$, 
one but not both of $\ext{N}/e$ and $\ext{N}\setminus e$ 
will
be the $\ext{0}$ extensor 
if and only if $e$ is a loop or a coloop in the matroid of $\ext{N}$.
If $\ext{N}'=\ext{0}$ then $\ext{N}'^{\perp}=\ext{0}$ and 
$\ext{M}_E(\ext{N}')=\ext{0}$.
We therefore write property 2 without restricting $e$ to a non-separator.
\item Property 1 except for signs is immediate from direct sum of 
subspaces and their corresponding extensors.  Property 2 except for the
signs
follows immediately 
from the fact that minor $[I_{\Is}V_{\Vs}E]$ of matrix (\ref{matrixM})
equals a linear combination with coefficients $g_e$ and $r_e$
because 
the column $e$ belongs to this minor
no matter which $e\in E$ is specified 
for the identity.
\end{enumerate}

\begin{proof}
From the definition, $\ext{M}(\ext{N}_1\ext{N}_2)$ $=$
$\Is(\ext{N}_1\ext{N}_2)\;\Vs\big((\ext{N}_1\ext{N}_2)^\perp\big)$
which equals
\[
\epsilon(S_1)\epsilon(S_2)\epsilon(S_1S_2)
(-1)^{\rho\ext{N}_1^\perp\rho\ext{N}_2}\;
\Is(\ext{N}_1\ext{N}_2)\;\Vs(\ext{N}_1^\perp\ext{N}_2^\perp)
\]
by Theorem \ref{Identities}(\ref{PerpPerp}).
$\Is(\ext{N}_1\ext{N}_2)\;\Vs(\ext{N}_1^\perp\ext{N}_2^\perp)$ $=$
$\Is(\ext{N}_1)\Is(\ext{N}_2)\Vs(\ext{N}_1^\perp)\Vs(\ext{N}_2^\perp)$
which equals (by (\ref{Anticommutative}))
\[
(-1)^{\rho\ext{N}_1^\perp\rho\ext{N}_2}
\;\Is(\ext{N}_1)\;\Vs(\ext{N}_1^\perp)\;\Is(\ext{N}_2)\;\Vs(\ext{N}_2^\perp).
\]
Therefore $\ext{M}(\ext{N}_1\ext{N}_2)$ $=$
\begin{equation}
\label{EqInMainThmProof1}
\epsilon(S_1)\epsilon(S_2)\epsilon(S_1S_2)
\;\Is(\ext{N}_1)\;\Vs(\ext{N}_1^\perp)\;\Is(\ext{N}_2)\;\Vs(\ext{N}_2^\perp).
\end{equation}
Therefore, 
$\ext{M}(\ext{N}_1\ext{N}_2)/E_1E_2$ $=$
\[
\epsilon(S_1)\epsilon(S_2)\epsilon(S_1S_2)\;
\big((\ext{M}(\ext{N}_1)/E_1)\;\ext{E_1}\;
(\ext{M}(\ext{N}_2)/E_2)\;\ext{E_2}\big)/E_1E_2
\]
because the only \Plucker coordinates of the form
$\ext{M}(\ext{N}_i)[X_i]$ for $i=1$ or $2$ that 
contribute to (\ref{EqInMainThmProof1}) when it is contracted
by $E_1E_2$ satisfy $E_i\subseteq X_i$.  
The anticommutativity law 
(\ref{Anticommutative}) then
implies that $\ext{M}(\ext{N}_1\ext{N}_2)/E_1E_2$ $=$
\begin{equation*}
\begin{split}
\epsilon(S_1)\epsilon(S_2)\epsilon(S_1S_2)(-1)^{\Card{E_1}\Card{P_2}}
\;\big(\ext{M}_{E_1}(\ext{N}_1)\;
 \ext{M}_{E_2}(\ext{N}_2)\;\ext{E_1}\;\ext{E_2}\big)/E_1E_2\\
=\epsilon(S_1)\epsilon(S_2)\epsilon(S_1S_2)(-1)^{\Card{E_1}\Card{P_2}}
\;\ext{M}_{E_1}(\ext{N}_1)\;\ext{M}_{E_2}(\ext{N}_2).
\end{split}
\end{equation*}
Since the sequence orders of the $S_i$ are arbitrary, let $S_i=P_iE_i$ for $i=1$ 
and $2$.
According to equation (\ref{PermuteContraction}),
$\ext{M}_E(\ext{N}_1\ext{N}_2)$ $=$ $\ext{M}(\ext{N}_1\ext{N}_2)/E$ $=$
$\epsilon(\sigma)\ext{M}(\ext{N}_1\ext{N}_2)/E_1E_2$ where $(E_1E_2)_\sigma$
$=$ $E$.  Therefore, $\ext{M}_E(\ext{N}_1\ext{N}_2)$ $=$
$\pm\ext{M}_{E_1}(\ext{N}_1)\;\ext{M}_{E_2}(\ext{N}_2)$ 
with the sign equal to
\begin{equation*}
\begin{split}
\epsilon(\sigma)
\epsilon(P_1E_1)\epsilon(P_2E_2)\epsilon(P_1E_1P_2E_2)
(-1)^{\Card{E_1}\Card{P_2}}\\
=
\epsilon(\sigma)
\epsilon(P_1E_1)\epsilon(P_2E_2)\epsilon(P_1P_2E_1E_2)\\
=
\epsilon^2(\sigma)
\epsilon(P_1E_1)\epsilon(P_2E_2)\epsilon(P_1P_2(E_1E_2)_\sigma)\\
=
\epsilon(P_1E_1)\epsilon(P_2E_2)\epsilon(P_1P_2E),
\end{split}
\end{equation*}
which proves property 1 of the theorem.

Now for property 2.  
Let us apply Theorem \ref{ExtDC} to $\ext{N}$ and $\ext{N}^\perp$, and 
apply $\Is$ and $\Vs$ respectively.
\[
\ext{N}=(\ext{N}/e)\ext{e} + (\ext{N}\setminus e)\ext{1}(e).
\]
\[
\ext{N}^\perp=(\ext{N}^\perp/e)\ext{e} + (\ext{N}^\perp\setminus e)\ext{1}(e).
\]
\begin{equation*}
\begin{split}
\Is(\ext{N})=&
\Is\big((\ext{N}/e)\ext{e}\big) + \Is\big((\ext{N}\setminus e)\ext{1}(e)\big)\\
            =&\Is(\ext{N}/e)g_e\ext{e} + \Is(\ext{N}\setminus e)\ext{1}(e).
\end{split}
\end{equation*}
\begin{equation*}
\begin{split}
\Vs(\ext{N}^\perp)=&\Vs
\big((\ext{N}^\perp/e)\ext{e}\big) + 
\Vs\big((\ext{N}^\perp\setminus e)\ext{1}(e)\big)\\
                  =&\Vs(\ext{N}^\perp/e)r_e\ext{e} + \Vs(\ext{N}^\perp\setminus e)\ext{1}(e).
\end{split}
\end{equation*}
The exterior product of $\Is(\ext{N})$ and $\Vs(\ext{N}^\perp)$ is therefore
\begin{equation*}
\begin{split}
&g_e\Is(\ext{N}/e)\;\ext{e}\;\Vs(\ext{N}^\perp\setminus e)+
r_e\Is(\ext{N}\setminus e)\;\Vs(\ext{N}^\perp/e)\;\ext{e}\\
+&g_er_e\;\Is(\ext{N}/e)\;\ext{e}\;\Vs(\ext{N}^\perp/e)\;\ext{e}\\
+&\Is(\ext{N}\setminus e)\;\Vs(\ext{N}^\perp\setminus e).
\end{split}
\end{equation*}
$\ext{M}_E(\ext{N})$ is the result of contracting the above by $E$.
The third term above is $\ext{0}$ because $\ext{e}$ is a repeated factor.
The last term will vanish when contracted by $E$ because none of 
its non-zero \Plucker coordinates have an index that contains $e$.  So we
will omit them in the following.
By Theorem \ref{Identities}(\ref{ContractPerp})
\[
\Vs(\ext{N}^\perp\setminus e) = \Vs\big((\ext{N}/e)^\perp\big)\;
    \epsilon(S')\epsilon(S'e)(-1)^{1\cdot(\Card{S}-\rho\ext{N})}
\]
and
\[
\Vs(\ext{N}^\perp/e)=\Vs\big((\ext{N}\setminus e\big)^\perp)\;
\epsilon(S')\epsilon(S'e).
\]
Notice that $\rho(\ext{N}/e^\perp)$ $=$ $\Card{S'}-(\rho\ext{N}-1)$ 
$=$ $\Card{S}-\rho\ext{N}$.  So $\ext{e}\;\Vs(\ext{N}/e)^\perp$ $=$
$\Vs(\ext{N}/e)^\perp\;\ext{e}(-1)^{\Card{S}-\rho\ext{N}}$.  Therefore when the
above substitutions are made we get
\[
\ext{M}_E(\ext{N})\epsilon(S')\epsilon(S'e) =
\big(g_e \Is(\ext{N}/e)\;\Vs((\ext{N}/e)^\perp)+
r_e\Is(\ext{N}\setminus e)\;\Vs((\ext{N}\setminus e)^\perp)\big)\;\ext{e}/E.
\]
Lemma \ref{lemmaEpsilonFacts} used with $\sigma$ such that
$E_\sigma$ $=$ $E'e$ shows that the right hand side is
\[
\epsilon(\sigma)((\cdots)\ext{e})/E'e = 
\epsilon(E)\epsilon(E'e)((\cdots)\ext{e})/E'e.
\]
So the right hand side is
\[
\epsilon(E)\epsilon(E'e)\big(g_e\ext{M}_{E'}(\ext{N}/e)+r_e\ext{M}_{E'}(\ext{N}\setminus e)\big).
\]
Since the sequence order of $S'$ is arbitrary, we can choose $S'=PE'$.  The sign correction is then
\begin{equation*}
\begin{split}
\epsilon(S')\epsilon(S'e)\epsilon(E)\epsilon(E'e)=\\
\epsilon(PE')\epsilon(PE'e)\epsilon(E)\epsilon(E'e).
\end{split}
\end{equation*}
Applying the permutation $\tau$ for which $(E'e)_\tau$ $=$ $E$ to the two
appearances of subsequence $E'e$ does not change this expression's value. 
Hence the sign correction is
\[
\epsilon(PE')\epsilon(PE)\epsilon(E)\epsilon(E)=\epsilon(PE')\epsilon(PE)
\]
and property 2 of the theorem is verified.

The definition of $\ext{M}_E$ immediately gives property 4, and, together
with the definition of extensor dual, gives property 3.
\end{proof}

\begin{corollary}
The set of $\ext{M}_E(\ext{N})$ 
obtained as the $g_e$, $r_e$ range over $\Reals$ for
each $e\in E$ represents the points in a 
projective subspace of a Grassmannian 
(which consists of all the
linear subspaces over $\Reals (P_{\Is}\dunion P_{\Vs})$ 
with 
dimension $\Card{P}$).
\end{corollary}

\begin{proof}
Induction: Use Theorem \ref{TheoremMain} property 2 for when $|E|>0$ and 
property 3
for the the basis.
\end{proof}

\begin{proposition}Given $\ext{N}=\ext{N}(P,E)$,
and sequences $I\subseteq P$, $V\subseteq P$, and $\scomp{V}=P\setminus V$,
\label{PropositionMain}
\[
\epsilon(\scomp{V}\scma V)\epsilon(PE)\ext{M}_E(\ext{N})[I_{\Is}V_{\Vs}] = 
\epsilon(P)\sum_{A\subseteq E}\ext{N}[IA]\ext{N}
[\scomp{V}A]g_Ar_{\scomp{A}}.
\]
\end{proposition}
\Remark: The only non-zero terms in this sum are those for which both $A\dunion I$
and $A\dunion \scomp{V}$ are bases in the matroid of $\ext{N}$.

\begin{proof}
Recalling that $E$ symbolizes a sequence $e_1\cdots e_n$, 
let $E_{i}=e_iE_{i+1}$ so $E_1=E$ and $E_{n+1}=\emptyset$.
When property (2) of Theorem \ref{TheoremMain}
is applied successively for $e:=e_i$, $E:=E_i$ and
$E':=E_{i+1}$ for $i=1,2,\ldots,n$
and the products are expanded, the result is a sum
of $2^n$ terms, one for each subset $A\subseteq E$.
For each fixed $i$, $1\leq i\leq n$,
the instances of 
the symbols $E$ and $E'=E\setminus e_i$ within 
all applications of property (2) each denote the same sequences.
Therefore, sign cancellation occurs and we can write
\begin{equation}
\label{ExpansionFormula}
\ext{M}_E(\ext{N})[I_{\Is}V_{\Vs}] = \epsilon(PE)
\sum_{A\subseteq E}\epsilon(P)g_Ar_{\bar{A}}
\ext{M}_{\emptyset}(\ext{N}/A\setminus\bar{A})[I_{\Is}V_{\Vs}].
\end{equation}
Property (3) combined with the definitions of extensor deletion,
contraction and dualization demonstrate that within each term
\[
\ext{M}_{\emptyset}(\ext{N}/A\setminus\bar{A})[I_{\Is}V_{\Vs}]
=
\epsilon(\bar{V}V)\ext{N}[IA]\ext{N}[\bar{V}A],
\]
and the conclusion follows.
\end{proof}

\begin{corollary}
\label{MENisNonZero}
If $\ext{N}\neq\ext{0}$ and the parameters are generic or are all
positive, then $\ext{M}_E(\ext{N})\neq\ext{0}$.
\end{corollary}
\begin{proof}
$\ext{N}$ has some basis, i.e., $B\subseteq P\dunion E$ for which $\ext{N}[B]\neq 0$.
(N.B. $B=\emptyset$ is possible.)   Take $I=B\cap P$ and $V=P\setminus I$.
Proposition \ref{PropositionMain} indicates 
$\ext{M}_E(\ext{N})[I_\Is V_\Vs]$ $\neq$ $\ext{0}$ 
with $A$ $=$ $B\setminus I$.
\end{proof}

\begin{corollary}
\label{CorollaryHomoDegree}
\Plucker coordinate $\ext{M}_E(\ext{N})[I_{\Is}V_{\Vs}]$ is a homogeneous
polynomial in the $g_e$, $r_e$ whose terms are square-free and have
degree $\Rank{\ext{N}}-\Card{I}$ in the $g_e$ and degree 
$\Card{E}-\Rank{\ext{N}}+\Card{I}$ $=$ 
$\Card{E}+\Card{P}-\Rank{\ext{N}}-\Card{V}$
$=$ $\Rank{\ext{N}^{\perp}}-\Card{V}$ in the $r_e$.
\end{corollary}

\begin{proof}
Immediate from Proposition
\ref{PropositionMain}, matroid duality and the fact $\ext{N}[X]\neq 0$
only if $\Card{X}=\Rank{\ext{N}}$.
\end{proof}

\begin{corollary}
\label{corollaryMESubsetSum}
\begin{equation}
\label{MESubsetSum}
\epsilon(PE)\ext{M}_E(\ext{N})=\epsilon(P)\sum_
             {\begin{array}{c}
                A\subseteq E:
	\rank_{\ext{N}}{A}=\Card{A},\\
	\rank{\ext{N}}-\rank{(\ext{N}/A|P)}-
        \rank_{\ext{N}}{A}=0
               \end{array}}
    \ext{M}_\emptyset(\ext{N}/A|P) g_Ar_{\scomp{A}}.
\end{equation}
\end{corollary}

\begin{proof}

The definition of extensor deletion indicates 
$\ext{N}/A|P$ is an alternative notation
for $\ext{N}/A\setminus\bar{A}$ when $A\subseteq E$,
$\bar{A}$ means $E\setminus A$ and the
ground set of $\ext{N}$ is $P\dunion E$.  The conditions
stated in Corollary \ref{CorollaryHomoDegree}
allow us to restrict the sum as indicated.
Hence
formula (\ref{ExpansionFormula}) for each \Plucker
coordinate is equivalent to the given expression 
for the extensor.
\end{proof}

The following definition and consequence of 
Proposition \ref{PropositionMain} clarify some 
of the sign behavior resulting from the definitions.

\begin{definition}
A function $F=F^\epsilon(X)$ whose value might depend on
the ground set orientation $\epsilon$ and on the sequence $X$ is said to be
\begin{enumerate}
\item \textbf{alternating in $X$} if $F^\epsilon(X_\sigma) = 
\epsilon(\sigma)F^\epsilon(X)$, for all
$\sigma\in\Perms_{\Card{X}}$; and
\item \textbf{alternating in $\epsilon$} if $F^{-\epsilon}(X)$ $=$
$-F^\epsilon(X)$.
\end{enumerate}
\end{definition}

\begin{corollary}
\label{BasisCountCorollary}
Let $Q\subseteq P_{\Is}\dunion P_{\Vs}$ with $\Card{Q}=\Card{P}$.
\begin{enumerate}
\item $\ext{M}^\epsilon_E(\pm \ext{N})[Q]$ is constant under sign change of
$\pm \ext{N}$, and is alternating in $E$, $\epsilon$ and $Q$.
\item $\epsilon(PE)\ext{M}^\epsilon_E(\pm \ext{N})[Q]$ 
is constant under sign
change of $\pm \ext{N}$ and under changes or reorderings of 
$\epsilon$ or  $E$; it
is alternating in $P$ and in $Q$.
\item $\epsilon(PE)\ext{M}^\epsilon_E(\pm \ext{N})[P_{\Is}]$
enumerates the bases of
$\mathcal{N}(\ext{N}/P)$, assuming $P$ is independent in the matroid 
$\mathcal{N}(\ext{N})$, by
\[
\epsilon(PE)\ext{M}^\epsilon_E(\pm \ext{N})[P_{\Is}]=
\sum_{B\subseteq E} g_Br_{\scomp{B}}\ext{N}^2[BP],
\]
\item and $\epsilon(PE)\ext{M}^\epsilon_E(\pm \ext{N})[P_{\Vs}]$ 
enumerates the bases of
$\mathcal{N}(\ext{N}\setminus P)$,  
assuming $P$ is coindependent in 
$\mathcal{N}(\ext{N})$, by\[
\epsilon(PE)\ext{M}^\epsilon_E(\pm \ext{N})[P_{\Vs}]=
\sum_{B\subseteq E} g_Br_{\scomp{B}}\ext{N}^2[B].
\]
\end{enumerate}
\end{corollary}

\Remark: Properties 3. and 4. express the Matrix Tree Theorem.

\section{Corank-Nullity Polynomials}
\label{Classify}

The well-known corank-nullity (or rank) polynomial is easily
generalized to ported matroids.  
We did\cite{sdcPorted} this
with a definition that differs from 
Las Vergnas' 
big Tutte polynomial \cite{SetPointedLV}
only in notation and in our applications.
We will now generalize it by including parameters and modify it
for oriented matroids by reinterpreting the symbols
for minors.
In the definition below
the symbol $\mathcal{N}/A|P$ represents the oriented minor
of the oriented matroid $\mathcal{N}$ obtained by contracting
$A$ and restricting to $P$.  Minors that are
different orientations of the same unoriented matroid are deemed
different objects.

\begin{definition}[Parametrized and Ported Corank-Nullity Polynomial]
\label{definitionRankPoly}
\[
R(\mathcal{N}(P,E))=\sum_{A\subseteq E}
	\MVAR{\mathcal{N}/A|P}g_Ar_{\scomp{A}}\;
	u^{\Rank{\mathcal{N}}-\Rank{[\mathcal{N}/A|P]}-\Rank{A}}\;
	v^{\Card{A}-\Rank{A}}.
\]
\end{definition}

In this formula, the bracketed oriented matroid 
$\MVAR{\mathcal{N}/A|P}$ $=$ 
$\MVAR{\mathcal{N}_{i1}\oplus\ldots\oplus\mathcal{N}_{ic}}$ denotes the
(commutative) product of the \textit{variables}
$\MVAR{\mathcal{N}_{i1}}$,  \ldots, $\MVAR{\mathcal{N}_{ic}}$, where each 
variable signifies a connected component of $\mathcal{N}/A|P$.  If
$P=\emptyset$ then $\MVAR{\mathcal{N}/A|P}$ $=$ 
$\MVAR{\emptyset}$ $= 1$; so
$R$ reduces to the corank-nullity polynomial, parametrized.

The formula 
therefore defines a polynomial in parameters $g_e,r_e$, whose 
(other) variables are $u$, $v$ together with a distinct variable for every
connected component of every minor of $\mathcal{N}$ obtained by
contracting some subset $A\subseteq E$ and deleting
$\scomp{A}=E\setminus A$.  The latter variables only occur
in monomials that signify direct sums of one or more minors.

It is readily verified that $R(\mathcal{N}(P,E))$ satisfies the 
ported Tutte equations below. 
The details
published in \cite{sdcPorted} can be immediately adapted 
to the changes we described above.
We state these results without proof:

\begin{proposition}
\label{RankPolyPortedTutteEq}
\begin{enumerate}
\item 
If $e\in E$ is neither a port nor a loop nor a coloop in $\mathcal{N}(P,E)$,
\begin{equation}
R(\mathcal{N}(P,E)) = g_eR(\mathcal{N}/e) + r_eR(\mathcal{N}\setminus e).
\end{equation}
\item
$R(\mathcal{N}_1\oplus\mathcal{N}_2) = R(\mathcal{N}_1)R(\mathcal{N}_2) =
R(\mathcal{N}_2)R(\mathcal{N}_1)$.
\item $R(\mathcal{N}_1(e))$ $=$ $g_e + r_eu$ and
$R(\mathcal{N}_0(e))$ $=$ $r_e + g_ev$, for the coloop and loop matroids
$\mathcal{N}_1(e)$ and $\mathcal{N}_0(e)$ on $E=\{e\}$, $P=\emptyset$.
\item
$R(\mathcal{N}(P,\emptyset)) = [\mathcal{N}]$ (i.e. when $E=\emptyset$.)
\end{enumerate}
\end{proposition}

Let us take $\mathcal{N}$ to be the oriented matroid
presented by extensor $\ext{N}$.
The reader can now verify that the expansion
in corollary \ref{corollaryMESubsetSum} 
for $\ext{M}_E(\ext{N})$
is obtained from
$R(\mathcal{N}(P,E))$ by substituting $u=0$, $v=0$ and the
extensor $\epsilon(P)\epsilon(PE)\ext{M}_\emptyset(\ext{N}/A|P)$
for the monomial $[\mathcal{N}/A|P]$ in the term with
factor $g_Ar_{\scomp{A}}$, for each $A$.  Note that 
$\Card{A}-\Rank{A}=0$ and 
$\Rank{\mathcal{N}}-\Rank{[\mathcal{N}/A|P]}-\Rank{A}$ 
$=\Rank{\mathcal{N}}-\Rank{(P\dunion A)} = 0$
imply that $A$ is independent and $P\dunion A$ is spanning
in $\mathcal{N}$. 
Therefore $\ext{N}/A|P\neq\ext{0}$ for those terms where
the exponents of $u$ and $v$ are both zero.

With arbitrary $\ext{N}$, the substitution
of $\epsilon(PE)\epsilon(P)\ext{M}_{\emptyset}(\ext{N}/A_i|P)$ for 
monomial $[\mathcal{N}/A_j|P]$ 
$=$ $[\mathcal{N}/A_i|P]$ 
in $R(\mathcal{N}(P,E))$ is not 
well-defined.  
The reason is that the same oriented matroid
$\mathcal{N}/A_j|P$ might be represented by different extensors all
with the form $\ext{N}/A_i|P$, for various $A_i\neq A_j$.  
They may differ by representing different subspaces 
(in the same oriented matroid stratification\cite[\textsection 2.4]{OMBOOK}
layer) of  $\FieldK P$.  
Therefore, different values $\ext{M}_E(\ext{N}/A_i|P)$ must be
substituted in $[\mathcal{N}/A_i|P]g_{A_i}r_{\scomp{A_i}}$ with
different $A_i$ even though these $[\mathcal{N}/A_i|P]$ all denote the same
oriented matroid.

The one general situation where $R(\mathcal{N})$ with $u=v=0$
determines $\ext{M}_E(\ext{N})$ is when $\ext{N}(P,E)$ is 
a unimodular extensor, i.e., one that represents the unimodular oriented 
matroid $\mathcal{N}$.  
\cite[Theorem 3.1.1, p. 41]{CombinatorialGeometries} 
provides this among other equivalent characterizations
of unimodular (also called regular) matroids.  
One of these characterizations is that
bracket values from $\{+1, -1, 0\}$ may be assigned so the 
Grassmann-\Plucker relationships hold over $\mathbb{Q}$.

\begin{definition}
\label{DefMUOM}
The extensor-valued function 
$\mathcal{N}(P,E)\rightarrow\ext{M}_E(\mathcal{N})$ is defined on the 
minor-closed class of ported unimodular oriented matroids by
\[
\ext{M}_E(\mathcal{N})=\ext{M}_E(\ext{N})
\]
where $\pm\ext{N}(P,E)$ are the two unimodular presentations of 
the ported oriented matroid $\mathcal{N}(P,E)$.
\end{definition}

The simplest case to demonstrate that the unimodular matroids
must be oriented in order for monomial substitution in
$R\big(\mathcal{N}(P,E)\big)$ to produce $\ext{M}_E(\mathcal{N})$ is
the two orientations $\mathcal{N}_1$, $\mathcal{N}_2$
of the 2-circuit matroid on two ports.  Here, $E=\emptyset$
and $\ext{M}_{\emptyset}(\mathcal{N}_1)$ $\neq$ 
$\ext{M}_{\emptyset}(\mathcal{N}_2)$.

We conclude:

\begin{theorem}
\label{MEofUnimodularExt}
For unimodular oriented matroid
$\mathcal{N}=\mathcal{N}(P,E)$, $\ext{M}_E(\mathcal{N})$ is the result of evaluating 
$R(\mathcal{N})$ (in the exterior algebra) after the substitutions
$u=0$, $v=0$ and extensor 
$\epsilon(P)\epsilon(PE)\ext{M}_\emptyset\big(\mathcal{N}_i(P)\big)$
for each monomial $[\mathcal{N}_i(P)]$ (which symbolizes an oriented 
unimodular matroid with ground set $P$).
\end{theorem}
\begin{proof}
Immediate from the above remarks and Corollary \ref{corollaryMESubsetSum}.
\end{proof}

\begin{theorem}
\label{EquationsExtFunMatroids}
The function $\ext{M}_E(\mathcal{N})$ defined above 
on ported unimodular oriented matroids satisfies the 
properties:  (Symbols like $E$ and $P$ 
denote sequences and $\ext{M}_E(\mathcal{N})$ 
depends on the $g_e, r_e$ and $\epsilon$.)
\begin{enumerate}
\item If $e\in E$ is neither a separator nor a port, and $E'=E\setminus E$,  
then
\[
\ext{M}_E(\mathcal{N})=\epsilon(PE)\epsilon(PE')\;
\big(g_e\ext{M}_{E'}(\mathcal{N}/e) + r_e \ext{M}_{E'}(\mathcal{N}\setminus e)\big).
\]
\item If 
$\mathcal{N}_1(P_1,E_1)$ and $\mathcal{N}_2(P_2,E_2)$ have 
disjoint ground sets and $E=E_1\dunion E_2$, then
\[
\ext{M}_E(\mathcal{N}_1\oplus\mathcal{N}_2)=
\epsilon(P_1P_2E)\epsilon(P_1E_1)\epsilon(P_2E_2)\;
\ext{M}_{E_1}(\mathcal{N}_1)\;\ext{M}_{E_2}(\mathcal{N}_2).
\]
\item If $P=\emptyset$ and  
$\mathcal{B}(\mathcal{N})$ denotes the collection of bases of 
$\mathcal{N}$, then
\[
\ext{M}_E(\mathcal{N})=\epsilon(E)\sum_{B\in\mathcal{B}(\mathcal{N})}g_Br_{\scomp{B}}.
\]

\end{enumerate}
\end{theorem}
\begin{proof}
Immediate from Theorem \ref{TheoremMain} and the above
remarks.
\end{proof}

\section{Basis, Set and Flat Expansions}
\label{Activity}

Theorem \ref{EquationsExtFunMatroids} shows 
that when $\mathcal{N}(P,E)$ is a unimodular matroid,
$\ext{M}_E(\mathcal{N})$ 
is a substitution of extensors
and $u=v=0$ into
$R_P(\mathcal{N})$ of Definition \ref{definitionRankPoly}.
Proposition \ref{RankPolyPortedTutteEq} demonstrates $R_P$
is a ported Tutte function.
This motivates general study of ported and parametrized 
Tutte functions of matroids.  Below, we extend some
known general expressions for Tutte functions and invariants.
Our results are
more clearly expressed and no harder to prove than if 
the parameters were omitted.  Furthermore,
we identify an unsolved problem due to 
combining parameters with ports.

Evidently, each Tutte function value $F(\mathcal{N})$
is determined via the Tutte
equations from parameter values and from values of $F$ on 
decomposable minors of $\mathcal{N}$.  One way to 
study a Tutte function is to present the definition independently 
of the Tutte equations and then prove that the function so 
defined satisfies the Tutte equations.  A second way is to
specify the parameters and the function values on indecomposable minors and
then prove a solution exists and is unique.  In other words, one
proves that all the Tutte equations are consistent with the given 
parameters and values, and that the solution is unique.  
(This extends to our Tutte functions remarks of Pak\cite{Pak}.)
The issue 
that arbitrary parameters and values on indecomposable sometimes
fail to be consistent or sometimes result in multiple solutions was studied
in \cite{MR93a:05047,BollobasRiordanTuttePolyColored,Ellis-Monaghan-Traldi}.
Fortunately, 
the Tutte functions (of matroids, oriented matroids and of extensors)
studied in this paper are all defined the first way.  However,
our work suggests the open problem to extend these studies to the ported
Tutte equations.  The indecomposables will then include connected matroids or 
oriented matroids with more than one element, i.e., not just
loops and coloops.

Zaslavsky\cite{MR93a:05047}
defines the normal class of parametrized (but not ported) Tutte
functions as those Tutte functions for which there exist 
$u$ and $v$ for which, for all $e$,
the point value on coloop $e$ is $r_eu + g_e$ 
and the point value
on loop $e$ is $g_ev + r_e$.  The normal Tutte functions 
are exactly those obtained by substitutions into the parametrized corank-nullity
polynomial.  All Tutte invariants are normal Tutte functions
and non-normal Tutte functions do not express much of
the matroid structure---See 
\cite{MR93a:05047,BollobasRiordanTuttePolyColored,Ellis-Monaghan-Traldi}
for details about how parametrization complicates Tutte invariant
theory.  

Our Theorem \ref{MEofUnimodularExt} expresses how
$\ext{M}_E(\mathcal{N})$ 
fits into the natural ported generalization of the
normal class.  Our extensor Tutte function of ported oriented 
unimodular matroids and its invariant specialization is expressed
by a substitution into the ported corank-nullity polynomial of
oriented matroids.

While only the normal Tutte functions have corank-nullity
polynomial expressions, they all have 
basis expansion expressions\cite{MR93a:05047}.
In the rest of this section, we discuss these and other
expansion expressions for ported unoriented and oriented matroids.

The basis expansion originated by Tutte\cite{TutteDich}
for graphs and Crapo\cite{CrapoAct} for matroids depends on 
a particular but arbitrary ground set element order $O$.
Each basis determines a term from the internal
and external activities of elements with respect to that basis
according to the ordering $O$.
Our way to generalize is to 
restrict $O$ to orders in which
every port element is ordered before 
each $e\in E$.  (We use the convention that the deleted/contracted element
is the last, i.e., greatest element under order $O$ eligible for reduction.)

Gordan and McMahon define\cite{GordonMcMachonGreedoid} 
a ``computation tree''
to formalize the application of a subset of Tutte equations
to a matroid and some of its minors.
Each (Tutte) computation
tree determines
a polynomial in the parameters and point values.
Therefore, when
$\mathcal{N}$ is in the domain of a Tutte function, each
of these computation trees determine the same value.  
Computation
trees are a way to give
a basis expansion expression in terms of 
a more 
general definition 
of internal and external activities of elements
with respect to a basis.  The expansion is more general because it
is based on any Tutte equation computation rather than on an element order $O$.
We will extend to ported computation trees the
classification\cite{GordonMcMachonGreedoid} of elements as internally or externally, 
active or passive
with respect to each path down the tree.
In each case, the result is an interval
partition of the boolean subset lattice of $E$.  

\begin{definition}
Given $\mathcal{N}(P,E)$,  a
\textbf{$P$-subbasis} $F\in \mathcal{B}_P(\mathcal{N})$
is an independent set  with $F\subseteq E$
(so $F\cap P=\emptyset$) for which $F\dunion P$ is a spanning set
for $\mathcal{N}(P,E)$ (in other words, $F$ spans $\mathcal{N}/P$,
see \cite{SetPointedLV}.)
\end{definition}

\begin{proposition}
For every $P$-subbasis $F$ there exists an independent set $Q\subseteq P$
that extends $F$ to a basis $F\dunion Q\in \mathcal{B}(\mathcal{N})$.
Conversely, if $B\in\mathcal{B}(\mathcal{N})$ then $F=B\cap E=B\setminus P$
is a $P$-subbasis.
\end{proposition}

\begin{proof} Immediate. \end{proof}

\begin{definition}[Activities with respect to a $P$-subbasis and an element
ordering $O$]
Let ordering $O$ have every $p\in P$ before every
$e\in E$.  Let $F$ be a $P$-subbasis.  Let $B$ be any basis for 
$\mathcal{N}$ with $F\subseteq B$.
\begin{itemize}
\item Element $e\in F$
is internally active if $e$ is the least element
within its principal cocircuit with respect to $B$.  Thus, this principal
cocircuit contains no ports.  The reader can verify this definition is 
independent of the $B$ chosen to extend $F$.  Elements $e\in F$ that are
not internally active are called internally inactive.
\item Dually, element $e\in E$ with $e\not\in F$ is externally 
active if $e$ is the least element within its principal circuit with
respect to $B$.  Thus, each externally active element is spanned by 
$F$.  Elements $e\in E\setminus F$ that are not externally active
are called externally inactive.
\end{itemize}
\end{definition}

\begin{definition}[Computation Tree, following \cite{GordonMcMachonGreedoid}]
\label{CompTreeDef}
A ported (Tutte) computation tree for $\mathcal{N}(P,E)$ is a
binary tree whose root is labeled by $\mathcal{N}$ and which satisfies:
\begin{enumerate}
\item If $\mathcal{N}$ has non-separating elements not in $P$, then 
the root has two subtrees and there exists one such element $e$ for which 
one subtree is a computation tree
for $\mathcal{N}/e$ and the other subtree is a computation tree for 
$\mathcal{N}\setminus e$.

The branch to $\mathcal{N}/e$ is labeled with ``$e$ contracted'' and 
the other branch is labeled ``$e$ deleted''.
\item Otherwise (i.e., every element in $S(\mathcal{N})\setminus P)$
is separating) the root is a leaf.
\end{enumerate}
\end{definition}

An immediate consequence is
\begin{proposition}
Each leaf of a $P$-ported computation tree for $\mathcal{N}(P,E)$
is labeled by the direct sum of some minor of $\mathcal{N}$ on $P$ 
(oriented if $\mathcal{N}$ is oriented) 
summed with loop and/or coloop matroids with 
ground sets $\{e\}$ for various distinct $e\in E$ (possibly none).
\end{proposition}

\begin{definition}[Activities with respect to a leaf]
\label{ActivityTreeDef}
For a ported computation tree for $\mathcal{N}(P,E)$, 
a given leaf, and the path from the root to this leaf:
\begin{itemize}
\item Each $e\in E$ labeled ``contracted'' along this path
is called \textbf{internally passive}.
\item Each coloop $e\in E$ in the leaf's matroid is
called \textbf{internally active}.
\item Each $e\in E$ labeled ``deleted'' along this path
is called \textbf{externally passive}.
\item Each loop $e\in E$ in the leaf's matroid is
called \textbf{externally active}.
\end{itemize}
\end{definition}

\begin{proposition}
Given a leaf of a ported computation tree for $\mathcal{N}(P,E)$:
The set of internally active or internally passive elements 
constitute a 
$P$-subbasis of $\mathcal{N}$ which we say 
\textbf{belongs to the leaf}.  
Furthermore, every $P$-subbasis $F$ of $\mathcal{N}$ belongs to a unique leaf.
\end{proposition}

\begin{proof}
For the purpose of this proof, let us extend Definition \ref{ActivityTreeDef}
so that, given a computation tree with a given node $i$ 
labeled by matroid $\mathcal{N}_i$,
$e\in E$ is called internally passive when $e$ is labeled 
``contracted'' along the path from root $\mathcal{N}$ to
node $i$.  Let $IP_i$ denote the set of such internally passive 
elements.

It is easy to prove by induction on the length of the root to node $i$ path
that
(1) $IP_i\cup S(\mathcal{N}_i)$ spans $\mathcal{N}$ and 
(2) $IP_i$ is an independent set in $\mathcal{N}$.  The proof
of (1) uses the fact that elements labeled deleted are non-separators.  The
proof of (2) uses the fact that for each non-separator 
$f\in\mathcal{N}/IP_i$, $f\cup IP_i$ is independent in $\mathcal{N}$.

These properties applied to a leaf demonstrate the first conclusion,
since each $e\in E$ in the leaf's matroid must be a separator by Definition 
\ref{CompTreeDef}.

Given a $P$-subbasis $F$, we can find the unique leaf as follows: Beginning
at the root, descend the tree according to the rule: At each branch node,
descend along the edge labeled ``$e$-contracted'' if $e\in F$ and along
the edge labeled ``$e$-deleted'' otherwise (when $e\not\in F$).
(This algorithm also operates on arbitrary $F'\subseteq E$.)
\end{proof}

The above definitions and properties enable us to conclude:
\begin{proposition}
Given element ordering $O$ in which every $p\in P$ is ordered
before each $e\not\in P$, suppose we construct the unique $P$-ported
computation tree $\mathcal{T}$ in which the greatest non-separator $e\in E$ is
deleted and contracted in the matroid of each tree node.

The activity of each $e\in E$ relative to ordering $O$ and
$P$-subbasis $F\subseteq E$ is the same as the activity
of $e$ defined with respect to the leaf 
belonging to $F$ in $\mathcal{T}$.
\end{proposition}

\begin{definition}
\label{ActivitySymbolsDef}
Given a computation tree for 
(oriented) matroid $\mathcal{N}(P,E)$,
each $P$-subbasis $F\subseteq E$
is associated with the following subsets of non-port elements
defined according to Definition \ref{ActivityTreeDef}
from the unique leaf determined by the algorithm given above.
\begin{itemize}
\item $IA(F)\subseteq F$ denotes the set of internally active elements,
\item $IP(F)\subseteq F$ denotes the set of internally passive elements,
\item $EA(F)\subseteq E\setminus F$ 
denotes the set of externally active elements,
and 
\item $EP(F)\subseteq E\setminus F$ denotes the set of externally
passive elements.
\end{itemize}
\end{definition}

\begin{proposition}
\label{PartitionProposition}
Given a computation tree for
$\mathcal{N}(P,E)$, 
the boolean lattice of subsets of $E$
is partitioned by the collection of
intervals $[IP(F),F\cup EA(F)]$ (note $F\cup EA(F)=IP(F)\cup A(F)$)
determined from the collection
of $P$-subbases $F$, which correspond to the leaves.
\end{proposition}

\begin{proof}
Every subset $F'\subseteq E=S(\mathcal{N})\setminus P$ belongs to the
unique interval corresponding to the unique leaf found by the tree descending
algorithm given at the end of the previous proof.
\end{proof}

Dualizing, we obtain:

\begin{proposition}
\label{DualPartitionProposition}
Given a computation tree for
$\mathcal{N}(P,E)$, 
the boolean lattice of subsets of $E$
is also partitioned by the collection of
intervals $[EP(F),E\setminus F\cup IA(F)]$ 
(note $E\setminus F\cup IA(F)$ $=$ $EP(F)\cup A(F)$).
\end{proposition}
\begin{proof} The dual of the tree descending algorithm is to descend along
the edge labeled ``$e$-deleted'' if $e\in F$.
\end{proof}


The following generalizes the basis expansion expression given
in \cite{MR93a:05047} to ported (oriented) matroids, as well as 
Theorem 8.1 of \cite{SetPointedLV}.

\begin{definition}
\label{TuttePolyExpression}
Given parameters $g_e$, $r_e$, point values $x_e$, $y_e$, and 
(oriented) $\mathcal{N}(P,E)$
the Tutte polynomial expression
determined by the sets in Definition 
\ref{ActivitySymbolsDef} 
from a computation tree is 
equal to
\begin{equation}
\sum_{F\in \mathcal{B}_P}[\mathcal{N}/F|P]
\;x_{IA(F)}\;g_{IP(F)}\;y_{EA(F)}\;r_{EP(F)}.
\end{equation}
\end{definition}

Each Tutte polynomial expression is constructed by applying some of the
Tutte equations.  Therefore, if $\mathcal{N}(P,E)$ is in the domain of
Tutte function $f$, then $f(\mathcal{N})$ is given by any Tutte polynomial
expression with $f(\mathcal{N}/F|P)$ substituted for each oriented or 
unoriented matroid monomial $[\mathcal{N}/F|P]$.  (This generalizes the 
expression used in \cite{MR93a:05047} to define the Tutte polynomial under
the condition that all expansions yield the same expression.)

\subsection{Boolean Interval Expansion}

The following proposition 
expresses the ported corank-nullity polynomial
in terms of a $P$-subbasis expansion.  It is obtained by substituting
binomials $x_e=g_e+r_eu$, $y_e=r_e+g_ev$ and leaving
the matroid variables unchanged
in Definition
\ref{TuttePolyExpression}.  The different expansions from
different element orderings and Tutte computation trees all 
express the same polynomial because Proposition
\ref{RankPolyPortedTutteEq} demonstrates that $R_P$ is a ported
Tutte function and the values of $R_P$ on coloop, loop and
indecomposable matroids are readily verified to be given
by these substitutions.  

\begin{proposition}
The polynomial $R_P(\mathcal{N})$ is given
by the following activities and boolean interval expansion formula:
\begin{equation}
R_P(\mathcal{N})=
\sum_{F\in \mathcal{B}_P}[\mathcal{N}/F|P]
\Big(
\sum_{\substack{
       IP(F)\subseteq K \subseteq F\\
       EP(F)\subseteq L \subseteq E\setminus F
      }}
 g_{K\cup (E\setminus F\setminus L)}\;
 v^{\Card{E\setminus F\setminus L}}\;
 r_{L\cup (F\setminus K)}\;
 u^{\Card{F\setminus K}}\;\;
\Big)
\end{equation}
\end{proposition}

\begin{proof}
Let $A=K\cup (E\setminus F\setminus L)$ within the above expansion.
We can verify $\scomp{A}=E\setminus A=L\cup (F\setminus K)$.
For each $A\subseteq E$ a unique $P$-subbasis $F$, and two tree leaves
are determined, one by the tree descending algorithm and the other leaf
by the dual algorithm.  Thus $A$ and $\scomp{A}$ respectively belong
to intervals within the boolean lattice partitions of Propositions
\ref{PartitionProposition} and \ref{DualPartitionProposition}.  In particular, 
$A\in[IP(F),F\cup EA(F)]$ and 
$\scomp{A}\in[EP(F),E\setminus F\cup IA(F)]$.
Therefore the terms in the above sum are equal one by one to the terms in
the corank-nullity polynomial's subset expansion
(Definition \ref{definitionRankPoly}).
\end{proof}

For the purposes of this paper it was sufficient to recognize
that our extensor valued ported parametrized Tutte function of 
unimodular oriented matroids belongs to the natural generalization 
of Zaslavsky's normal class.  As such, it has, for arbitrary
parameters, expressions obtained by substitutions into
(1) computation trees, (2) ported parametrized Tutte polynomials
from such trees, (3) various $P$-subbasis expansions, 
and (4) the ported parametrized corank-nullity polynomial.

We state here the open problem to include ports into the results of 
Zaslavsky, Bollobas and Riordan, and Ellis-Monaghan and Traldi:
Can we classify with universal forms all of the ported parametrized Tutte
functions according to their parameters, non-port
point values, and the 
values on oriented or unoriented minors on port ground sets?  

\subsection{Geometric Lattice Flat Expansion}

A formula for the unparametrized ported Tutte (or corank-nullity) 
polynomials
of non-oriented matroids in terms of the lattice of flats (closed sets)
and its Mobius function was given in \cite{sdcPorted}.  
We generalize: (1) The expansion's
monomials $[\mathcal{Q}]$ can signify either oriented matroid minors,
when $\mathcal{N}$ is oriented, or non-oriented minors
when $\mathcal{N}$ is not oriented.  (2) The polynomial is parametrized
with $r_e$, $g_e$ for each $e\in E$.
The derivation relies on the
fact that the oriented or non-oriented matroid minor $[\mathcal{N}/A|P]$
(according to whether $\mathcal{N}$ is oriented or not) depends
only on the flat spanned by $A\subseteq E$.  

\begin{proposition}
Let $\mathcal{N}(P,E)$ be an oriented or unoriented.
Let $R_P(\mathcal{N})$ be given from Definition 
\ref{definitionRankPoly}.
In the formula below,
$F$ and $G$ range over the geometric lattice of flats contained
in $E$.
\begin{equation}
R_P(\mathcal{N})(u,v) = \sum_{\mathcal{Q}}[\mathcal{Q}]
      \sum_{\substack{F\leq E\\
                     [\mathcal{N}/F|P]=[\mathcal{Q}]
           }}
                   u^{\rank{\mathcal{N}}-\rank{\mathcal{Q}}-\rank{F}}
                   v^{-\rank{F}}
                   \sum_{G\le F}
                   \mu(G,F)
                   \prod_{e\in G}
                    (r_e+g_ev)
\end{equation}
\end{proposition}

\begin{proof}
It follows the steps for theorem 8 in \cite{sdcPorted}.
\end{proof}

\Remark: The chirotope values for the oriented matroid minor
$\mathcal{N}'$ $=$ $\mathcal{N}/ F|P$ are
$\chi_{\mathcal{N}'}(X)$
$=$ $\chi_{\mathcal{N}}(XB_F)$ where $X$ is restricted to sequences
over $P$ and 
$B_F$ is any basis for the flat spanned by $F$.
While this formula
defines a chirotope function 
only up to a constant sign factor, the oriented matroid (which is what
the monomial $[\mathcal{N}/F|P]$ denotes) is uniquely defined.  We mention
this because when we evaluate the corank-nullity polynomial to obtain
$\ext{M}_E(\mathcal{N})$ for the unimodular oriented matroid 
$\mathcal{N}$
substitute an extensor 
for each $[\mathcal{N}/A|P]$.  
However, the object we
substitute is $\ext{M}_\emptyset(\ext{N}/A|P)$,
not $\ext{N}/A|P$.  It is the unique extensor defined by equation 
\eqref{MENComponentDef} applied to one of the chirotopes
that present $[\mathcal{N}/A|P]$ (or to any other representation
of $[\mathcal{N}/A|P]$ for that matter).  We already remarked that equation 
\eqref{MENComponentDef} is unchanged
when its argument changes sign!

\section{Electrical Networks}
\label{ENets}
The most common way to formulate the electrical network problem 
that occurs in our subject is to use the discrete and parametrized
Laplace's equation.  The insights from matroid theory appear 
with greater clarity when the analyses of the solution, which is
in terms of the electrical potential values at graph vertices, are 
replaced by analyses in terms of voltage and current variables 
directly associated with edges.  After we review the Laplacian,
we move on to edge-based formulations.  Some details of relationships
between our results and electrical network problems are then demonstrated.

The following observations underlie the combination of ideas
about electrical networks and ported oriented matroids:
First, the basis enumerator is a Tutte function that 
happens to be a determinant in the case of unimodular 
(or regular) matroids, such as the graphic matroids.
Second, topics involving determinants, including the
chirotope presentation of realizable oriented matroids,
can profitably be studied with exterior algebra.  
Third,
non-trivial electrical networks (see \textsection \ref{ElecNetEquatSec}) 
must have port elements 
for supplying power or for connecting to an external environment,
in addition to the resistor elements usually modeled by
graph edges.  We are interested in their combinatorial 
properties beyond spanning tree 
counts\cite{sdcMTT,sdcBDIMatroid,sdcELEW,sdcPorted,sdcISCAS95,sdcOMP,sdcISCAS98}.
Finally, electrical flows and potential differences
are inherently directional.  The patterns of their directions
(expressed by sign functions on graph edges)
feasible under Kirchhoff's two laws are precisely the 
vector and covector families, respectively, of the graphic 
oriented matroid.  Indeed, the duality between these laws, current
and voltage, is characterized by oriented matroid theory.

\subsection{Discrete Laplacian}
\label{Laplacian}
The combinatorial (or discrete) Laplacian
is the matrix of coefficients
in the equations (\ref{LaplaceEquation}) 
below
in variables $\phi_i$, $1\leq i \leq n$.
These discrete Laplace equations model (among other 
situations) a resistive electrical network
when $\phi_i$ represents the electrical potential 
(or voltage) at vertex $i$ and
constant $I_i$ represents the current flowing into
vertex $i$.
\begin{equation}
\label{LaplaceEquation}
\sum_{\{j | e=ij\in E\}}
g_e(\phi_i-\phi_j) = I_i \;\;\; 1\le i \le n
\end{equation}

If each  conductance $g_{e}$ is non-negative, 
or is either zero
or generic, then the rank of 
the $n\times n$ Laplacian matrix is $n-k$, 
where $k$ is the number of connected components in the $n$
vertex undirected graph whose edges
are the $e=ij$ with $g_{e}\neq 0$.  Each order $n-k$ non-singular
diagonal submatrix is called a reduced Laplacian.  
(In \textsection \ref{ElecNetEquatSec}, we will express Laplace's equations
with an $\ext{M}_E(\ext{N})$ independent of $k$.) 
The reduced Laplace 
equations, together with $\phi_i=0$ for each vertex $i$ corresponding
to a deleted column, model a network where each such $i$ is 
contracted into a single 
\textbf{grounded} vertex
whose potential is fixed to zero and whose 
external current is 
unrestricted\footnote{Do not confuse with ground set.}.  
Since the graph is undirected, equations (\ref{LaplaceEquation})
imply that the current into the grounded vertex equals
the sum of the $I_j$ for the non-grounded vertices.

The inverse of a reduced Laplacian matrix is called the
discrete Green's function in \cite{ChungYaoGreens}.
This inverse matrix's elements are each expressed (using Cramer's rule)
by a ratio of an order $n-k-1$ minor to a common order $n-k$ minor
denominator.
The 
list of all minors, of all orders,
is an example of \Plucker coordinates---Here, these are the 
$\binom{2n}{n}$ maximal minors of the matrix obtained by appending
the $n\times n$ identity matrix to the side of the Laplacian.

The Matrix Tree Theorem asserts that each $n-1$ order minor 
equals $\pm\sum g_T$, the enumeration of spanning trees $T$ by products
of edge parameters $g_T=\prod_{e \in T}g_e$.  
See \cite{sdcMTT} for 
similar interpretations of all the minors and 
for generalizations to directed graphs.
The formulas we call ``Maxwell's rules'' were given without
proof by Maxwell\cite{MaxR}, and the dual forms of them were proved
by Kirchhoff\cite{Kirchhoff}.
Maxwell also described the static equilibrium solution for 
stressed linear elastic
framework in terms of enumerations over minimally rigid 
subframeworks\cite{MaxwellsFramesPaper}; this enumerated set is the basis
set for the rigidity matroid\cite{RigidityBook}.  The one-dimensional case
is analogous to the electrical problem.
The survey by Biggs \cite{BiggsAlgPotTheory}
covers the discrete Laplacian, the Matrix Tree Theorem, and the
use of spanning tree enumeration to solve the discrete Laplace equations, and
many additional topics, including the asymmetric discrete Laplacian.
Biggs presents the Kirchhoff's solution method.
and Nerode and Shank\cite{NerodeShank}, also used by 
Bott and Duffin\cite{BottDuffinAlgNetworks}, Smith\cite{SmithElec}
and Maurer\cite{Maurer76}.
This method constructs a symmetric projection matrix from a sum of fundamental
cocycle matrices, one for each spanning tree.  Analysis of basis exchange,
i.e., the pivot calculation
implies the appropriately weighted matrix sum is symmetric.  We plan to present
the generalization of this argument to extensors in a future publication.

Tree counting, the discrete Laplacian and electrical network models
with parameters 
have a spectrum of applications including  electrical circuit theory,
knot theory, random walks and the analysis of Markov chains
(see for example
\cite{DoyleSnellRandom,BiggsAlgPotTheory,ChungYaoCovHeatKer,ChungYaoGreens}).
Their application to square dissections
is described in \cite{BSST,TutteBook}; Tutte gives a Laplacian based 
``barycentric embedding'' proof
of Kuratowski's Theorem in \cite{HowToDraw}.

It is generally known among electrical engineers in circuit theory
that the same kinds of homogeneous rational polynomial functions that appear
in Maxwell's rule occur generally as the coefficients (and 
minors of them) in all of the 
linear relationships between the port quantities that define the
externally observable characteristics of a linear resistive network.
Our results display this principle within the mathematical contexts of 
the enumerative combinatorics of graphs, oriented matroids and 
exterior algebra: See Corollary \ref{CorollaryHomoDegree}.
Some electrical network analysis software actually enumerates trees and related
structures to do ``symbolic analysis.''
See for example 
\cite{TraWilBilinear,ThulasiramanSwamyBook,someEEBooks,sdcISCAS98}.  

\subsection{Electrical Network Equations in terms of Edges}
\label{ElecNetEquatSec}
Consider a graph with two kinds of edges, called ports $P$ and
resistors $E$.  The graph is directed with an arbitrary edge 
orientation.  Let unimodular extensor $\ext{N}(P,E)$ present
its ported graphic oriented matroid.  Let $r$ be the rank of
this matroid.

Let $g_e,r_e$ be parameters for each $e\in E$.  The extensors
$\Is(\ext{N})(P_\Is\dunion P_\Vs\dunion E)$
and 
$\Vs(\ext{N}^\perp)(P_\Is\dunion P_\Vs\dunion E)$ defined in \textsection
\ref{Extensors} determine the electrical network equations in the
way expressed by Theorem \ref{OrthogSpace} applied to $\ext{N}^\perp$.

To be specific, these equations are a linear system on the
$\Card{E}+2\Card{P}$ variables $\{x_e,\ldots;$
$ v_p,\ldots;$
$ i_p,\ldots\}$.
Each $e\in E$ is associated to variable $x_e$.  Each $p\in P$ is 
associated to two variables, $v_p$ called the voltage and 
$i_p$ called the current.
Let matrix $K$ with $r$ rows be any matrix that presents
$\Is(\ext{N})$.   Let $C$ be any matrix
with $\Card{P}+\Card{E}-r$ rows that presents $\Vs(\ext{N}^\perp)$.  
These matrices express 
Kirchhoff's equations combined with a homogeneous expression of Ohm's
law.
$K$ determines the following \textbf{current} equations:
\[
\sum_{p\in P}K_{j,p_\Is} i_p + \sum_{e\in E}K_{j,e}x_e = 0
\text{ for } j = 1, \ldots, r.
\]
$C$ determines the following \textbf{voltage} equations:
\[
\sum_{p\in P}C_{j,p_\Vs} v_p + \sum_{e\in E}C_{j,e}x_e = 0
\text{ for } j = 1, \ldots, \Card{E}+\Card{P} -r.
\]
It is helpful to see the electrical network equations in terms of 
$\ext{N}$ directly.  Let $N$ be any matrix presentation of $\ext{N}$; 
for example, a reduced oriented incidence matrix of the graph.
Let $N^\perp$ present $\ext{N}^\perp$; the rows of $N^\perp$ comprise a
basis for the cycle space of the graph.  The current equations can 
be written:
\[
\sum_{p\in P}N_{j,p} i_p + \sum_{e\in E}N_{j,e}g_ex_e = 0
\text{ for } j = 1, \ldots, r.
\]
The voltage equations can be written:
\[
\sum_{p\in P}N^\perp_{j,p} v_p + \sum_{e\in E}N^\perp_{j,e}r_ex_e = 0
\text{ for } j = 1, \ldots, \Card{E}+\Card{P} -r.
\]
The equations which $\ext{M}_E(\ext{N})$ 
presents are obtained by eliminating all 
the variables $x_e$, $e\in E$ in
the voltage and current equations taken together.  
Corollary \ref{MENisNonZero} indicates the rank of the resulting system
of $\Card{P}$ equations on $2\Card{P}$ variables is $\Card{P}$, provided that 
the parameters are generic or all positive.

The above analysis illustrates the role for the port element
distinction in modeling a physical system.
Each non-port 
element models a completely defined subsystem.
The ``proto-voltage\cite{SmithElec}'' $x_e$ parametrizes the
state of one electrical resistor, for example.
The behavior of this resistor is thus defined by 
Ohm's law:  
When $g_e$ and $r_e$ are both non-zero, the current is $g_ex_e$ if and only if the voltage is $r_ex_e$.
The entire model (the graph, for example) specifies all the interactions 
(via Kirchhoff's laws, for electricity) between
its subsystems.
Each port element models an 
interface pertaining to an interaction 
of the system with an unspecified environment, for observing the system
behavior of interest to the application, and to help specify how 
certain larger systems are composed out of 
previously entire 
subsystems.  
For us, the environment
is assumed, for each port,
to constrain the currents into one terminal and out of the other
terminal to be equal. Environmental constraints between
voltages at terminals belonging to distinct ports are forbidden as well.
(Engineering models encompass multiport elements, whose
behavior is specified using multiple port elements\cite{ChensBook,Recski,LinNonLinCircuitsBook}.
For example, a linear multiport element is specified a linear
constraints among the variables associated with its ports; this generalizes
Ohm's law to so-called multi-terminal resistors.  
Each of our
ported objects 
can model a single multiport element
within a larger model.
A topic for 
future research is to abstract this along the lines given here.)

Let a graph on vertices $\{1,\ldots, n\}$ be given with conductance parameters 
$g_e$ for each edge.  
We now derive Laplace's equations from the voltage and current equations.
Let us
append a new vertex $0$ (which will be grounded) and 
$n$ port edges $p_i\in P$, with
each $p_i$ directed from vertex $0$ to vertex $i$.  For simplicity take
parameter $r_e=1$ for each edge $e$ from the original graph.  We can choose
$N$ so that the current equations are 
\[
i_p = \sum_{e\in E}J_{p,e}g_ex_e,\;\;p\in P
\]
and the voltage equations are
\[
x_e = \phi_{p_r} - \phi_{p_s} = \sum_{i=1}^n J_{p_i,e} \phi_{p_i},
\text{ where }e=rs, e\in E.
\]
where we used potential (relative to vertex $0$) $\phi_i$ in place of 
$v_{p_i}$ and $J$ is the oriented vertex-edge incidence matrix of the original
graph.  Laplace's equation is obtained by eliminating 
the variables $x_e$ which represent differences of potential across
resistor edges.  Indeed, one presentation for 
$\ext{M}_E(\ext{N})$ in this case is the $n\times 2n$ matrix
$[I_n \Delta]$ formed by concatenating the identity matrix with
the Laplacian matrix $\Delta$.  One manifestation of Theorem 
\ref{TheoremMain} is therefore that each of the forest enumerating
polynomials given by an arbitrary minor, of any order, of the Laplacian
is a (non-ported) Tutte function of graphic oriented matroids.

\subsection{Maxwell's Rules}
\label{MaxwellSection}

Kirchhoff\cite{Kirchhoff} first described the solution of a resistive electrical
network problem in terms of enumerations of spanning trees
and certain forests (or dual forests).  One such description
is the following surprising yet
classical result, one of the ``easily remembered''
rules stated without proof by Maxwell\cite{MaxR}.

Let's discuss Maxwell's rule for a network with
one port edge $p$ in terms of graph edges.
Let all edges $\neq p$ represent unit resistors in an
electrical network.  Port edge $p$ is to 
demarks two terminal vertices.  Maxwell's rule asserts 
that the equivalent resistance between the two terminal vertices 
equals the quotient of the count $A$ of spanning trees that
contain edge $p$ divided by the count $B$ of spanning trees 
that omit edge $p$.  Note that whenever a two-tree spanning forest $F$ 
for which $F\dunion p$ is a tree counted by $A$ is contracted and the 
other non-$p$ edges deleted,
$p$ becomes a coloop.  When this is done for $B$, $p$ becomes
a loop.  $A$ and $B$ are 
the two \Plucker coordinates of the extensor we construct in this 
paper.  They are the coefficients in the linear equation 
$Ai_p + B(-v_p) = 0$ that relates
the port current and the (negated) port voltage.  
The reader can verify that each of $A$ and $B$ satisfy 
the additive ported Tutte equation;
so $(A,B)$ satisfies it.  
(It is instructive to calculate how $(A,B)$ behaves under the 
multiplicative ported Tutte equation.)  The spanning tree count
resistance expression $(A,B)$ generalizes with weighted spanning 
tree generating function when conductance parameters $g_e$ are given
for $e\neq p$.

Suppose  two ports, $p_1$ and $p_2$ are given.  
One, say $p_1$, denotes the terminal pair
between which a specified current value $i_{p_1}$ is constrained to
flow (assuming this is feasible.)  The other, $p_2$, denotes the terminals
between which to observe or measure the potential difference 
$v_{p_2}$ (voltage)
that results when the specified current constraint is the only 
electrical power source (assuming this voltage value is unique.)
Maxwell's rule for one port is the specialization with $p_1=p_2=p$.

The two-port Maxwell's rule holds that the port variables are related
by $Ai_{p_1} + B(-v_{p_2})=0$ were $A$ and $B$ are again generating
functions, except that some of the monomials in $A$ might be negative.
The spanning trees that omit $p_1$ and $p_2$ are enumerated by $B$.
Term $A$ enumerates forests $F$ for which each of $F\cup p_1$ and 
$F\cup p_2$ is a spanning tree are enumerated by $A$.  For each such
$F$, $F\cup\{p_1, p_2\}$ contains a unique circuit, which of course
contains both $p_1$ and $p_2$.  The relative 
directions in which $p_1$ and $p_2$ are traversed determines the
sign of $\pm g_F$ in $A$.

It is amusing go back to the one port $p$ case and 
derive 
the Maxwell's rule coefficients $(A,B)$ from the principles we just
illustrated.
The port behavior of the network consisting of one coloop $p$ is
defined by $1\cdot i_p + 0\cdot v_p =0$: This constraint 
can be expressed by 
\Plucker coordinates
$(1,0)$.  The port current is 0 and the port voltage is
unconstrained since the graph has no cycles.  
Dually, the port behavior of the one port loop network
is expressed by $(0,1)$.  The port voltage is 0 but the current is 
unconstrained.  One consequence of our main theorem proves that
$A(1,0)+B(0,1)=(A,B)$ gives the  
\Plucker coordinates of the constraint that the original graph
imposes on its port current and (negated) port voltage variables,
in accordance with Maxwell's rule $Ai_p + B(-v_p) = 0$.
Our theory justifies the choices of $(1,0)$ and $(0,1)$
rather than any other multiple $(\alpha,0)$ or $(0,\beta)$ in
$A(\alpha,0)+B(0,\beta)$.

\subsection{Deriving Maxwell's Rule}

See \cite{ChensBook}
for an elementary derivation that includes the
version of Maxwell's rule\cite{MaxR} that applies to 
to a pair of
ports that do not share a common vertex.  
Graphic matroid orientation becomes
relevant in this situation.
Explicit port edges
have proven their usefulness in electrical network 
analysis\cite[\textsection 13.6]{LinNonLinCircuitsBook}.
In this situation,  some
\Plucker coordinates of 
$\ext{M}_E(\ext{N})$, as polynomials in $g_e,r_e$, have
\emph{terms of opposite sign} only when two port edges do not share a 
vertex.  Our contribution is to characterize the signs within the theory of 
oriented matroids and Tutte functions.
It is true that such polynomials can be expressed in terms
of minors of the Laplacian; this was done by manipulation
of solutions to Laplace's equation in \cite{BSST}; see also \cite{TutteBook}.
However, our extensor and oriented matroid formulation enables 
the analysis to be done without the introduction of vertices.

Our derivation of Maxwell's rule for two ports (of which the
one port version is a special case) begins with the electrical network
equations with $P=\{p_1,p_2\}$.  Let $\ext{M}=\ext{M}_E(\ext{N})$ as defined
in \textsection \ref{Extensors} be as discussed above, and let $M$ be
any $2\times 4$ matrix presentation of $\ext{M}$.  The
two equations
\[
M\left[\begin{array}{c}i_1 \\  i_2  \\  v_1 \\  v_2\end{array}\right]
= 0
\]
are obtained by eliminating the variables $x_e$, $e\in E$ from the
electrical network equations.  The currents $i_1$, $i_2$ in edges
$p_1=ab$ and $p_2=cd$ flow from vertices $a$ to $b$, and $c$ to $d$
respectively.  The voltage (drop) $v_1$ across edge $p_1$ is the 
potential at $a$ minus the potential at $b$; the corresponding
convention defines the voltage $v_2$ across edge $p_2$.  

We will assume that all the $r_e=1$ and 
that $\ext{M}[p_{1\Vs}p_{2\Vs}]\neq 0$.
The latter is assured from Corollary \ref{BasisCountCorollary} (4.) 
provided that $E$ contains
a spanning tree and all the $g_e$ are either positive or generic.  Under
these conditions, the \textbf{transfer resistance} $\rho_{21}$ given
by $(-v_2)/i_1$ when $i_2=0$ and $i_1\neq 0$ is well-defined.  (These 
conventions are used so that 
when $p_1$ and $p_2$ are identical or parallel,
$\rho$ signifies the familiar equivalent resistance which is always positive
when $E$ is connected and all $g_e>0$.)

\begin{proposition}[Maxwell's Rule]
\label{Maxwell1}
Given the electrical network graph model described above, let
$\mathcal{B}$ denote the collection of edge sets 
$T\subseteq E$ of trees that span the vertex set $V$, and assume 
$\sum_{T\in \mathcal{B}}g_T\neq 0$.

For vertices $i,j,k,l$, 
let $\mathcal{B}_{ik,jl}$ be the collection 
of all $F\subseteq E$
for which
the subgraph $(V,F)$ is a forest
with exactly two trees where vertices
$i$ and $k$ are in one tree and 
$j$ and $l$ are in the other tree.

The transfer resistance $\rho_{21}$, where $p_1=ab$ and 
$p_2=cd$, is well-defined and is given by:

\begin{equation}
\label{Maxwell1equation2port}
\rho_{21}=\frac{\sum_{F\in \mathcal{B}_{ac,bd}}g_F -
                  \sum_{F\in \mathcal{B}_{ad,bc}}g_F}
                 {\sum_{T\in \mathcal{B}}g_T}.
\end{equation}
\end{proposition}

\begin{proof}
We will abuse the notation slightly by using $v_k$ and $i_k$ for
ground set elements $p_{k\Vs}$ and $p_{k\Is}$, and $\ext{v}_k$,
$\ext{i}_k$ for the corresponding extensors, $k=1,2$.

Corollary \ref{BasisCountCorollary} (4.) shows that 
$\ext{M}[v_1v_2]$ ($=$ 
$\ext{M}[p_{1\Vs}p_{2\Vs}]$) $=$ 
$\sum_{T\in \mathcal{B}}g_T$ $\neq$ $0$, so by Cramer's rule,
\[
\rho_{21}=
-\frac{v_2}{i_1}=-\Big(-\frac{\ext{M}_E(\ext{N})[v_1i_1]}{\ext{M}_E(\ext{N})[v_1v_2]}\Big).
\]
Let us apply Corollary \ref{corollaryMESubsetSum} to the numerator and 
denominator.  To do this, we first calculate 
$\ext{M}_{\emptyset}(\ext{N})$ for 6 unimodular extensors
$\ext{N}(\{p_1,p_2\},\emptyset)$ that present the 6 oriented matroids
on $P=\{p_1,p_2\}$, which are all graphic.  For each of the 6,
we can then determine
the \Plucker coordinate 
values $\ext{M}_{\emptyset}(\ext{N})[v_1i_1]$
and $\ext{M}_{\emptyset}(\ext{N})[v_1v_2]$.  Those oriented matroid
minors for which one of these values is non-zero
will characterize, together with the rank conditions in
Corollary \ref{corollaryMESubsetSum}, which forests or trees contribute to
each sum.  These characterizations of the forest or tree terms, and
of their signs, analyzed for each
\Plucker coordinate, will complete the derivation of Maxwell's rule.

\def\omN{\mathcal{N}}

Four of the 6 oriented matroids are the direct sums
of either the loop $\omN_0(p_1)$ or coloop $\omN_1(p_1)$ 
with either the loop or coloop on $p_2$.   
The other two oriented matroids are the orientations
of the 2-circuit matroid on $\{p_1,p_2\}$.
Let $\omN_1^+$ denote the oriented circuit
$\pm(+-)$; graphically, $p_1$ and $p_2$ are parallel.
So $\omN_1^-$ denotes the oriented matroid of antiparallel
$p_1$ and $p_2$.
Table \ref{6Table} lists 
these six distinct 
ported oriented matroids $\mathcal{N}(\{p_1,p_2\},\emptyset)$,
their unimodular extensor presentations $\ext{N}(\{p_1,p_2\})$,
and the corresponding extensor values $\ext{M}_{\emptyset}(\ext{N})$.
\begin{table}[!ht]
\caption{\label{6Table}The six oriented matroids on $\{p_1,p_2\}$.}
\[
\begin{array}{c|c|c}
\mathrm{matroid} & \ext{N} & \ext{M}_\emptyset(\ext{N}) \\ \hline
\omN_0(p_1)\oplus \omN_0(p_2) &\pm 1& \ext{v}_1\ext{v}_2 \\
\omN_0(p_1)\oplus \omN_1(p_2) &\pm \ext{p}_2& \ext{i}_2\ext{v}_1 \\
\omN_1(p_1)\oplus \omN_0(p_2) &\pm \ext{p}_1& \ext{i}_1\ext{v}_2 \\
\omN_1(p_1)\oplus \omN_1(p_2) &\pm \ext{p}_1\ext{p}_2& \ext{i}_1\ext{i}_2 \\
\omN_1^- &\pm(\ext{p}_1-\ext{p}_2)& (\ext{i}_1-\ext{i}_2)(\ext{v}_1+\ext{v}_2)=
    \ext{i}_1\ext{v}_1+\ext{i}_1\ext{v}_2-\ext{i}_2\ext{v}_1-\ext{i}_2\ext{v}_2 \\
\omN_1^+ &\pm(\ext{p}_1+\ext{p}_2)& (\ext{i}_1+\ext{i}_2)(\ext{v}_2-\ext{v}_1)=
     \ext{i}_1\ext{v}_2-\ext{i}_1\ext{v}_1+\ext{i}_2\ext{v}_2-\ext{i}_2\ext{v}_1
\end{array}
\]
\end{table}
The $\ext{M}_{\emptyset}(\ext{N})$ 
values are easily found up to sign.  The signs are given in
Proposition \ref{PropositionMain}.  
One way to calculate is to analyze 
the corresponding electrical network with
2 ports and no resistors.
For example, the network with
oriented matroid $\omN_0(p_1)\oplus \omN_1(p_2)$ 
constrains its port $p_1$, a loop, to
have  voltage drop $v_1=0$ but its port $p_2$, a coloop, to have current
$i_2=0$.
The current in the loop and voltage across the coloop are unconstrained.
The solution subspace 
corresponds to equations 
($v_1=0$; $i_2=0$).
The extensors representing these equations
are 
$\alpha \ext{v}_1\ext{i}_2$, $\alpha\neq 0$.

Similarly, the network of two parallel ports (case $\omN_1^+$)
constrains the sum of
voltage drops going around the oriented circuit to be $0$, so Kirchhoff's
voltage law is expressed by $v_1-v_2=0$.  Kirchhoff's current law
in the same network is expressed $i_1+i_2=0$.  Hence the corresponding
extensor is $\pm(\ext{v}_1-\ext{v}_2)(\ext{i}_1+\ext{i}_2)$.

We complete the derivation.
First for the denominator.
From table \ref{6Table} the only terms in 
(\ref{MESubsetSum}) of Corollary \ref{corollaryMESubsetSum} 
that might contribute to 
$\ext{M}_E(\ext{N})[v_1v_2]$ are those for which $\omN/A|P$ is the matroid
of two loops $\omN_0(p_1)\oplus \omN_0(p_2)$ because the only appearance
of $\ext{v}_1\ext{v}_2$ is in that matroid's row.  The rank conditions further
restrict the contributing $A$ to spanning trees.

Finally, for the numerator $\ext{M}_E(\ext{N})[v_1i_1]$, we locate 
$\pm \ext{v}_1\ext{i}_1$ in the bottom two rows.
These appearances have \textit{opposite} sign.
For $A\subseteq E$ with $\omN/A|P=\omN_1^-$, the contribution
is $-g_A$.  
The sign is opposite
when $\omN/A|P=\omN_1^+$, so the distinct orientations of the 2-circuit
obtained when contracting $A$ account for the opposite signs
in (\ref{Maxwell1equation2port}).  We can again verify from the rank
conditions that the $F$ $=$ $A$ contributing to the numerator of 
(\ref{Maxwell1equation2port}) are the spanning forests with 2 trees
containing the indicted vertices as claimed.

Note that the sign dependence of $\ext{M}_E(\ext{N})$ on $\epsilon$ and
the order of $P=p_1p_2$ cancels in the ratio $\rho_{21}$.
\end{proof}

While Theorem \ref{Maxwell1} can be proved by elementary arguments as in
\cite{ChensBook}, the above proof demonstrates how it can be derived
from the forgoing theory using algebraic calculations.

\Remark: The one port version 
is immediately derived using a graph where 
$p_1$ and $p_2$ are parallel edges because our proof puts no
special conditions on the two ports.

\subsection{Signed Contributions of Spanning Forests}

In some cases of graphs given with two or more port edges, some of these
\Plucker coordinates will equal the difference 
between the counts of two kinds of spanning forests.  Such a coordinate
pertains to the coefficient that relates a quantity observed at 
one port to of a voltage or a current quantity at a different port.  
Our results show how oriented matroid properties determine the sign by which
each forest $F$ contributes to this coordinate 
(Corollary \ref{corollaryMESubsetSum}).
In particular,
that sign is determined by the graphic oriented matroid on the (directed)
port edges obtained by contracting $F$ and deleting the remaining non-port
edges.  
The simplest case where distinct signs do occur is when both orientations
of the same 2-circuit matroid on two ports appear in this process.
The contribution of $F$ to the \Plucker coordinate with 
a given index $X$ is 
calculated in a particularly simple way:  We solve the electrical network 
\textit{of port edges only} (with no resistors!)
that resulted from this deletion and contraction,
after checking matroid rank conditions necessary for at least
one coordinate to be non-zero.
The contribution \emph{equals} the \Plucker coordinate with the same index 
$X$ from
the solution of the latter, port-edge-only network, weighted by
$g_Fr_{E\setminus F}$.
It should be noted that the sign of the contribution
is determined by the oriented graphic matroid on ports only, independent
of the particular $F$ contracted to obtain this oriented matroid.
Details appeared in our proof of Maxwell's rule for two ports, section 
\ref{MaxwellSection}.

\section{Additional Background and Directions}
\subsection{Ground Set Orientation}
\label{GroundSetOrientationBack}

The ground set orientation
and its role in defining a
canonical dual of an extensor, and our $\ext{M}_E(\ext{N})$ are motivated by
the idea of orientations of orientable manifolds and the definition
of pseudo-forms (or ``forms of odd-kind'' attributed to de Rham in
\cite{Frankel}) in the mathematics of physics.  A pseudo-form
is an antisymmetric multilinear operator $f=f_\epsilon$
that is parametrized by the orientation $\epsilon$ and for which
$\epsilon f_\epsilon$ is independent of the orientation\cite{Frankel}.
So, $\epsilon f_\epsilon$ is a well-defined 
form.  In physics, an orientation
specifies one's convention, say by a right-handed coordinate system, 
for how one defines
a positive volume or other naturally unsigned physical quantity
in terms of an exterior algebra form.  The orientation specifies
which ordered bases determine right handed coordinate systems.
In this context, the orientation $\epsilon$ is a $\pm 1$ function 
for which 
$\epsilon(B_1)\epsilon(B_2)$
is the sign of the determinant of the local Jacobian 
matrix which relates the ordered bases $B_1,B_2$.  

\subsection{Computational Complexity}
\label{Complexity}

Among the non-trivial Tutte
invariant functions of 
succinctly presented graphs or matroids,
only two (unless $\mathcal{P}=\mathcal{\#P}$)
are polynomial time 
computable
\cite{JaegerVertiganWelsh,MR99k:05053}.
One such function,
the number of bases,
is computable by the Matrix Tree Theorem
for graphs and its extension 
to unimodular matroids.  
This number 
is well-known as the evaluation $T(\mathcal{N},1,1)$ of the 
Tutte polynomial
function $T(\mathcal{N},x,y)$ of matroids $\mathcal{N}$.
However, 
computing $T(\mathcal{N},1,1)$ 
is $\#\mathcal{P}$-complete 
for arbitrary non-unimodular matroids\cite{VertiganBaseCounting}; this follows
because counting the perfect matchings in a bipartite graph is 
a $\#\mathcal{P}$-complete problem\cite{ValiantNumPComplete}.

The other easy-to-compute invariant is
determined by the 
dimension of the intersection of a linear subspace and its orthogonal
complement over a finite field\cite{MR99k:05053}.
The cited papers prove
that all of the other Tutte
matroid invariants are
either trivial or $\#\mathcal{P}$-complete.
More recently, analogous computational hardness results have been
proven for
Tutte functions of graphs (thus implying 
their hardness for matroids).
Among these results, 
is that  evaluating
the parametrized Tutte polynomial for given matroids
is a \textbf{VNP}-complete problem\cite{MakowskyLotzComplexity}.
Here, Valiant's non-uniform algebraic complexity model\cite{ValiantNonUnif}
is used, which 
counts as one 
deterministic
step each evaluation of a
polynomial on
constants, variables or previously
computed values.
\textbf{VNP} is this model's class that is analogous to 
\textbf{NP} in the Turing machine model.
(See the references in \cite{MakowskyLotzComplexity}).

A full account of 
the computational 
complexity
of Tutte
invariants of graphs and matroids is given in
\cite{JaegerVertiganWelsh,MR94m:57027,MR99k:05053,OxleyWelshTuttePolyPolyTime,
WelshRandomApprox,WelshApproxCounting}.

We remark for computationally-inclined readers that: 
\begin{enumerate}
\item The Tutte equations
describe non-unique recursive algorithms to compute Tutte functions
that generally require $2^\Card{E}$
steps. 
\item  A $\Card{P}\times 2\Card{P}$ matrix representing our
extensor 
can be computed from a 
graph or locally or totally unimodular 
matrix presentation of a ported oriented unimodular matroid.

One suitable algorithm is simple matrix block manipulations followed by
Gaussian elimination.  Such elimination-based algorithms
use polynomial bounded numbers of field operations.  Therefore, computation
of our extensor generalization of the basis enumerator on graphic
and other unimodular matroids is a polynomial time
problem when all $r_e=g_e=1$.

\end{enumerate}

\subsection{Other Directions}
\label{Peripheral}

It is natural to generalize the current and voltage equations 
so their respective solutions subspaces (taken to be within $\FieldK S$)
are 
not orthogonal\cite{sdcBDIMatroid}.  
This leads to the directed graph version of the Matrix Tree
Theorem.  It did lead as well to a ``oriented matroid pair'' model for 
combinatorial conditions
for certain equations with monotone non-linearities to be uniquely 
solvable\cite{sdcOMP}.  These conditions were stated in terms of
two oriented matroids with a common ground set
having complementary rank and no common non-zero covector; the current paper 
provides the insight that these two were obtained by deletion/contractions
to eliminate port elements.  Investigations of a generalization of the 
Tutte polynomial to two matroids with a common ground set were also
begun in \cite{WelshKayibiLinking}.

The computation tree formalism was used in
greedoid generalizations\cite{GordonMcMachonGreedoid}.
of the Tutte polynomial because those generalizations do not
always have an activities expansions based on element orders.
We leave investigation of ``ported greedoids'' to the future.

\section{Acknowledgments}

The author wishes to thank Prof. Ariel Caticha, Physics Department of
the University at Albany, for his invitation to participate in his
graduate Geometry of Physics course; Joanna Ellis-Monaghan, Michel Las
Vergnas, Gary Gordon, Alan D. Sokal, Richard P. Stanley, Lorenzo
Traldi and David G. Wagner for helpful discussions, correspondence and
encouragement; and Thomas Zaslavsky for his assistance leading to a
much improved presentation of the subject.

\bibliographystyle{abbrv}
\bibliography{TutteEx}

\end{document}